\documentclass[12pt]{article}
\usepackage{amsmath,amsfonts,amssymb,amsthm,mathrsfs}
\usepackage[a4paper,vmargin={3.5cm,3.5cm},hmargin={2.5cm,2.5cm}]{geometry}
\usepackage[font=sf, labelfont={sf,bf}, margin=1cm]{caption}
\usepackage{graphicx,graphics}
\usepackage{epsfig}
\usepackage{latexsym}
\usepackage[applemac]{inputenc}
\usepackage{ae,aecompl}
\usepackage[english]{babel}
 \usepackage[colorlinks=true]{hyperref}
\usepackage{pstricks}
\usepackage{enumerate}

\date{}

\newtheorem{theorem}{Theorem}[]

\newtheorem{proposition}[theorem]{Proposition}
\newtheorem{lemma}[theorem]{Lemma}

\def\llbracket{[\hspace{-.10em} [ }
\def\rrbracket{ ] \hspace{-.10em}]}

\title{\bf \textsc{The stable trees are nested}}

\author{\text{Nicolas Curien}\thanks{ \'Ecole normale supérieure, E-mail: nicolas.curien@ens.fr}  \ \ \& \text{B\'en\'edicte  Haas}\thanks{ Universit\'e Paris-Dauphine and \'Ecole normale supérieure, E-mail: haas@ceremade.dauphine.fr} }

\begin{document}
\maketitle
\abstract{We show that we can construct simultaneously all the stable trees as a nested family. More precisely, if $1 < \alpha < \alpha' \leq 2$ we prove that hidden inside any $\alpha$-stable tree we can find a  version of an $\alpha'$-stable tree rescaled by an independent Mittag-Leffler type distribution. This tree can be explicitly constructed  by a pruning procedure of the underlying stable tree or by a modification of the fragmentation associated with it. 
Our proofs are based on a recursive construction due to Marchal which is proved to converge almost surely towards a stable tree.}

\medskip


\medskip


\begin{figure}[!h]
 \begin{center}
 \includegraphics[width=16cm]{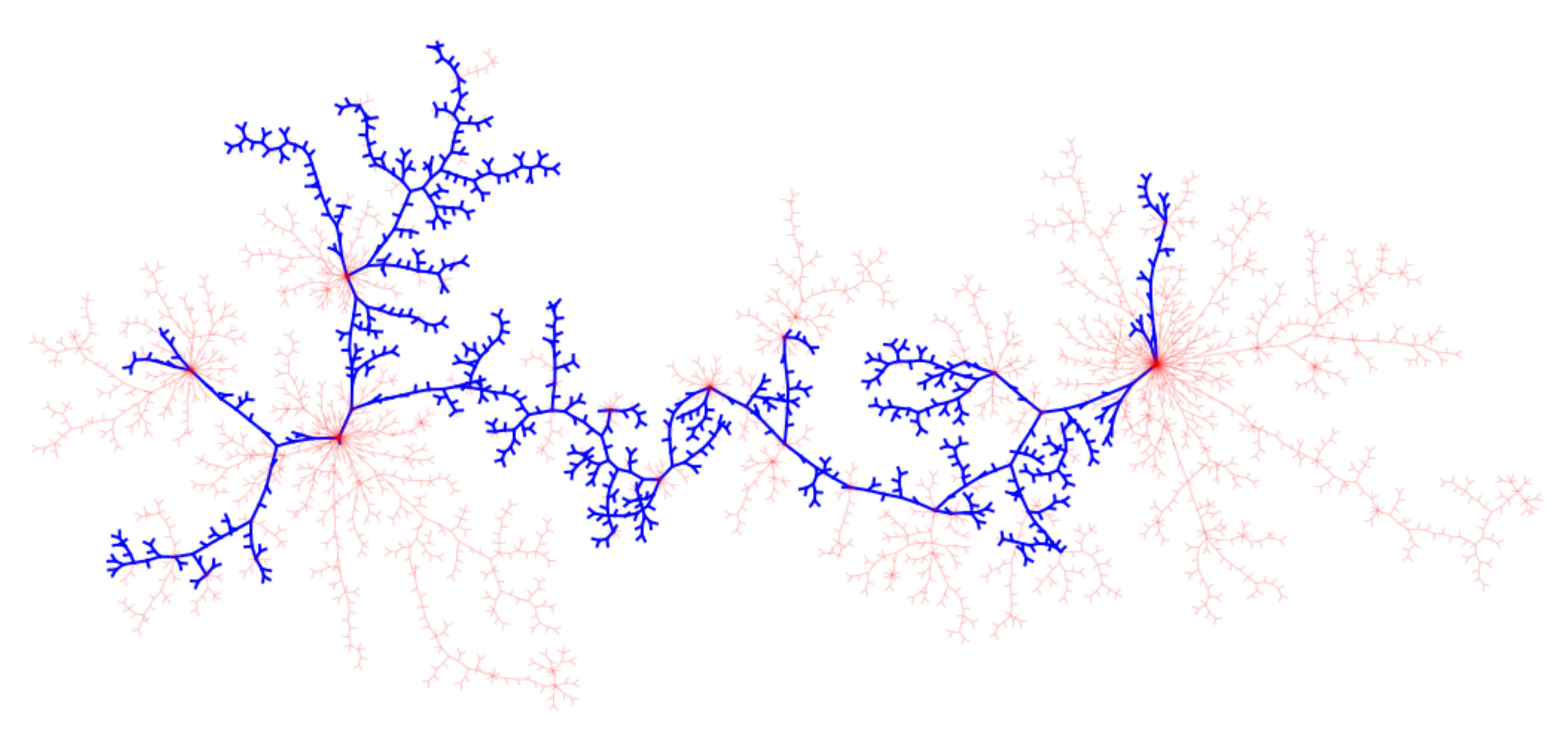}
 \caption{ \label{fig:inside}a Brownian tree hidden in a stable tree}
 \end{center}
 \end{figure}

\clearpage

\section*{Introduction}
Since the beginning of the early 1990's and the introduction of the Brownian continuum random tree by Aldous \cite{Ald91a}, random trees have been the object of an intense research in probability theory. Various classes of continuous random trees have been considered to extend the initial Brownian tree setup, two important such classes are Lévy and fragmentation trees. The Lévy trees have been introduced by Le Gall and Le Jan  \cite{LGLJ98} to describe the genealogical structure of continuous-state branching processes and furnish all the possible scaling limits of Galton-Watson trees \cite{DLG02}.  The fragmentation trees \cite{HM04} encode the genealogy of self-similar fragmentation processes and are the scaling limits of Markov branching trees \cite{HM12}. The intersection of these two classes consists of the stable trees of parameter
 \begin{eqnarray*} \alpha &\in& (1,2].  \end{eqnarray*}
 When $\alpha=2$, the $2$-stable tree corresponds to Aldous' Brownian CRT, up to a multiplicative scaling. For general $\alpha \in (1,2]$, the $\alpha$-stable trees can either be seen as the scaling limits of conditioned critical Galton-Watson trees with offspring distribution in the domain of attraction of a stable law of index $\alpha$ \cite{Du03}, or as a $1/\alpha-1$ self-similar fragmentation tree which is invariant under uniform re-rooting \cite{HPW09}. The increasing function $\alpha \mapsto 1-1/\alpha$ will repeatedly appear when dealing with stable trees and
we shall use the following notation throughout the paper: 
 \begin{eqnarray*}\bar{\alpha} &=& 1-1/\alpha. \end{eqnarray*}
We adopt the convention of Duquesne and  Le Gall \cite{DLG02} and denote by $ \mathcal{T}_{\alpha}$ the  \textit{standard} $\alpha$-stable tree, which is the one describing the genealogical structure of a continuous-state branching process with branching mechanism $\lambda \mapsto  \lambda^{\alpha}$. Alternatively, for $\alpha=2$, it can be defined as  the scaling limit of a conditioned Galton-Watson tree with offspring distribution $\eta$ with mean 1 and  variance $\sigma^2 \in (0,\infty)$  in the following sense:
 \begin{eqnarray*}
 \frac{T_n^{\mathrm{GW}}}{\sqrt n} &\xrightarrow[{n \to \infty}]{(\mathrm d)}& \frac{\sqrt 2}{\sigma} \mathcal T_2,
  \end{eqnarray*}
where $T_n^{\mathrm{GW}}$ is a version of the above mentioned Galton-Watson tree conditioned to have $n$ vertices. Similarly, for $\alpha \in (1,2)$ and for an offspring distribution $\eta$ with mean 1 and
such that $\eta(k)\sim Ck^{-1-\alpha}$ as $k \rightarrow \infty$, we have that:
 \begin{eqnarray*}
 \frac{T_n^{\mathrm{GW}}}{n^{\bar{\alpha}}} &\xrightarrow[{n \to \infty}]{(\mathrm d)}& \left(\frac{\alpha(\alpha-1)}{C\Gamma(2-\alpha)}\right)^{1/\alpha} \mathcal T_{\alpha}.
 \end{eqnarray*}
 These convergences hold for the Gromov-Hausdorff topology, as recalled in Section \ref{SecStable}.

\bigskip

\noindent \textbf{The nested family of stable trees.} Some connections between different stable and Lévy trees have already been observed. In \cite{ADV10}, Abraham, Delmas and Voisin  present a  pruning procedure for a large class of Lévy trees that leads to other Lévy trees. However, when applied to stable trees, their procedure gives pruned subtrees that are not stable anymore. In the other direction,  Bertoin, Le Gall and Le Jan \cite{BLGLJ07} obtain stable trees by folding  randomly the branches of a Brownian tree. These results are a priori unrelated to the ones presented here.

It is well-known \cite{DLG05,HM04} that the Hausdorff dimension of the $\alpha$-stable tree is almost surely   $1/\bar \alpha$,  for all  $\alpha \in (1,2]$. Thus,  in some sense, the stable trees are decreasing in the parameter $\alpha$.   The goal of this work is to give a precise geometric statement of the last heuristic and to show that all the stable trees can be seen as a single family of (random) nested trees.  To do so, rather than considering the standard versions of the stable trees $ \mathcal{T}_{\alpha}$, we will consider  randomly scaled versions  so that it is possible to build on a same probability space  a family of nested stable trees. More precisely, we let  $\Gamma_{a}$, $a>0$, denote a Gamma random variable of parameter $a$, that is with density proportional to $x^{a-1}e^{-x}\mathbf{1}_{\{x >0\}}$, and set  \begin{eqnarray*} J_{\alpha} & \overset{( \mathrm{d})}{=} & \alpha \cdot \left(\Gamma_{1+\bar \alpha}\right)^{\bar \alpha}.  \end{eqnarray*}
This variable has been designed so that if we rescale a stable tree $ \mathcal{T}_{\alpha}$ by a independent variable distributed as $ J_{\alpha}$ (that is we multiply the distances in $ \mathcal{T}_{\alpha}$ by the factor $ J_{\alpha}$) then the height of a random uniform point in $ J_{\alpha} \cdot \mathcal{T}_{\alpha}$ has a distribution that does not depend on $\alpha$, namely  a Gamma distribution of parameter $2$. Our main result is then

\begin{theorem}
\label{th:processus}
There exists a process of rescaled nested stable trees $(\mathscr{T}_{\alpha},1 < \alpha \leq 2)$ such that:
\begin{enumerate}
\item[$\bullet$] $\mathscr{T}_{\alpha}\overset{(\mathrm{d})}=J_{\alpha} \cdot \mathcal T_{\alpha}$, \ where $J_{\alpha}$  is independent of $\mathcal T_{\alpha}$,  $\alpha \in (1,2],$
\item[$\bullet$] $\mathscr{T}_{\alpha'} \subset \mathscr{T}_{\alpha}$, \ for all $1<\alpha \leq \alpha ' \leq 2$.
\end{enumerate}
\end{theorem}
The existence of this decreasing process of rescaled stable trees  is actually a simple corollary of the following proposition, which says that hidden inside any stable tree of parameter $\alpha$, there exists a rescaled version of a stable tree of parameter $\alpha'>\alpha$   (see Fig.\,\ref{fig:inside}).   
\begin{proposition} \label{thm:inside} Let $1 < \alpha < \alpha' \leq 2.$ There exists a closed (random) subtree $ \mathfrak{T}_{\alpha,\alpha'}$ of the $\alpha$-stable tree $\mathcal{T}_{\alpha}$,  such that   \begin{eqnarray*} \mathfrak{T}_{\alpha,\alpha'} & \overset{(\mathrm{d})}{=}& M_{\alpha,\alpha'}^{\bar \alpha'} \cdot \mathcal{T}_{\alpha'},  \end{eqnarray*} where $ \mathcal{T}_{\alpha'}$ is a standard $\alpha'$-stable tree and $ M_{\alpha,\alpha'}$ is an independent variable distributed as $(\alpha'/\alpha)^{1/\bar \alpha'}$ times a generalized Mittag-Leffler distribution  with parameters $\big({\bar \alpha}/{\bar \alpha'}, \bar \alpha \big)$. \end{proposition}
 
The definition of generalized Mittag-Leffler distributions is recalled in Section \ref{sec:RM}.  
 Thorough this work, the boundary case $\alpha' =2$, where a (rescaled) Brownian tree is extracted from an $\alpha$-stable tree, deserves a special attention since our constructions and proofs simplify substantially in this case. Although not unique,  the subtree $ \mathfrak{T}_{\alpha,\alpha'}$ of the last proposition can be explicitly  constructed by a pruning procedure of the tree $ \mathcal{T}_{\alpha}$.

\bigskip 

\begin{figure}[!h]
 \begin{center}
 \includegraphics[width=12cm]{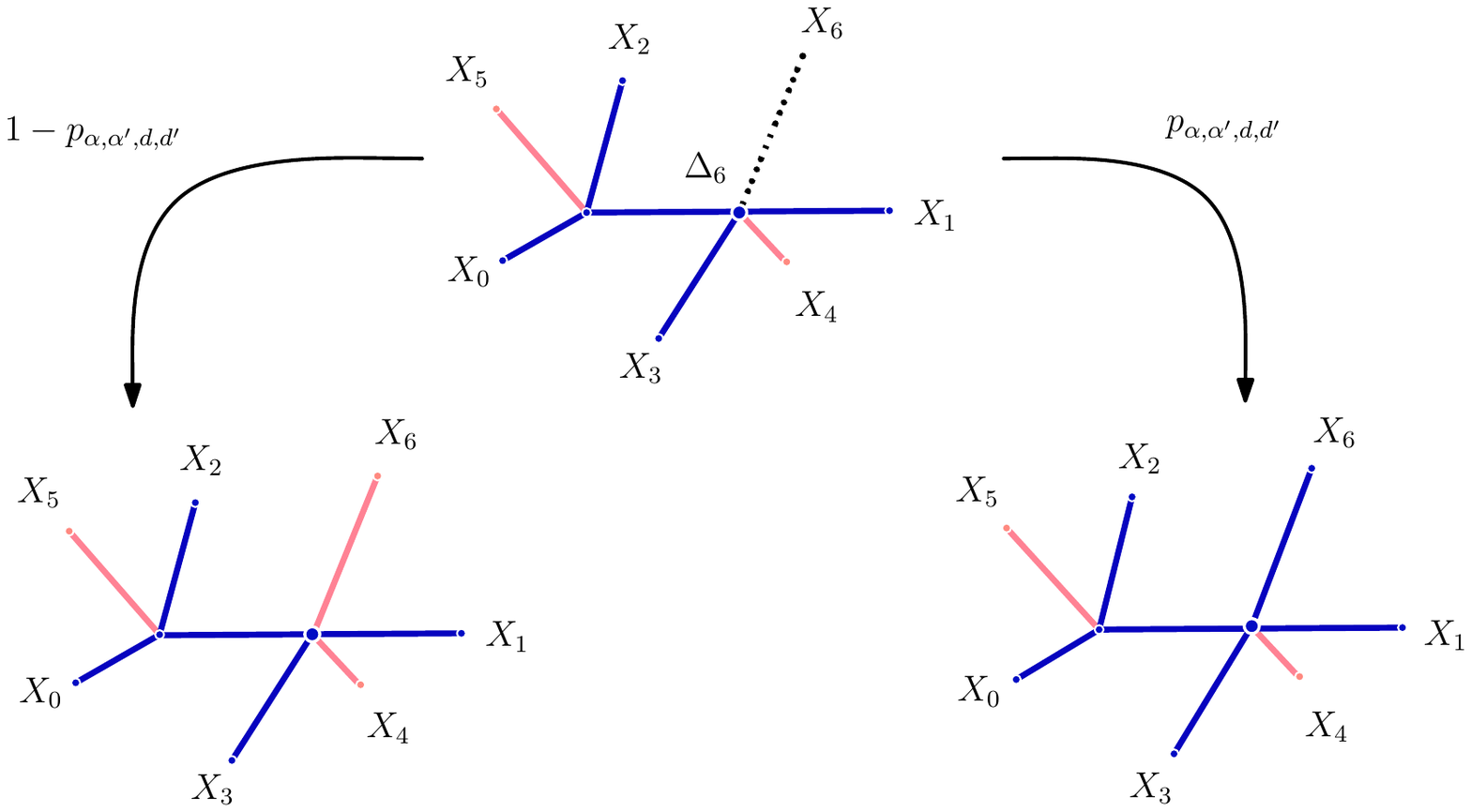}
 \caption{ \label{fig:binary} illustration of the pruning construction. Here $d'=3$, $d=4$ and the trees $\tau_5$ and $\tau_6$ are in blue.}
 \end{center}
 \end{figure}

\noindent \textbf{Pruning procedure.} Let  $1 < \alpha < \alpha' \leq 2.$
Conditionally on the stable tree $ \mathcal{T}_{\alpha}$, let $X_{i},i \geq 0$ be i.i.d. random leaves sampled according to its uniform mass measure. Then write $ \mathfrak{t}_n$ for the subtree spanned by $X_0, X_1, ..., X_n$ in $ \mathcal{T}_\alpha$.   For $d,d' \in \{0,2,3,4,5,...\}$ with $d \geq d'$ and $d \geq 2$, we introduce the following probabilities 
  \begin{eqnarray} \label{pdd} p_{\alpha,\alpha',d,d'} &=& \left\{ \begin{array}{ccl}
  0 & & \mbox{if }d'=0,\\
  1 && \mbox{if } d'=2,\\
\displaystyle \frac{(d'-1-\alpha')(\alpha-1)}{(d-1-\alpha)(\alpha'-1)}&& \mbox{otherwise}, \end{array} \right. \qquad \in [0,1], \end{eqnarray} from which we  construct inductively a sequence of subtrees $\tau_{n} \subset  \mathfrak{t}_n \subset \mathcal{T}_{\alpha}$ as follows: \begin{itemize}
\item $ \tau_1=  \mathfrak{t}_1 = \llbracket X_0, X_1 \rrbracket$  is the line segment linking $X_0$ to $X_1$ in $ \mathcal{T}_{\alpha}$,
\item  at step $i\geq 2$, write $ \mathfrak{t}_{i} = \mathfrak{t}_{i-1} \cup \llbracket X_{i}, \Delta_{i}\rrbracket$, where $\Delta_{i}\in \mathfrak{t}_{i-1} $ and $\llbracket X_i,\Delta_i \rrbracket$ is the shortest path in $\mathfrak{t}_i$ connecting $X_{i}$ to $ \mathfrak{t}_{i-1}$. Then let $d_{i}$ denote the degree (multiplicity) of $\Delta_{i}$ in $ \mathfrak{t}_{i-1}$ and $d'_{i}$  its degree in $ \tau_{i-1}$, with the convention that $ d_i'=0$ if $\Delta_{i} \notin \tau_{i-1}$. Finally set  (see Fig.\,\ref{fig:binary})  \begin{eqnarray*}  \tau_{i} = \tau_{i-1}\cup \llbracket X_{i}, \Delta_{i} \rrbracket   & \mbox{ with probability }& p_{\alpha,\alpha',d_{i},d'_{i}}, \\
 \tau_{i}= \tau_{i-1} &\mbox{otherwise}.  \end{eqnarray*}
\end{itemize}
 The sequence $(\tau_n)$ is clearly increasing in $ \mathcal{T}_{\alpha}$ and we denote by $ \mathrm{Prun}_{\alpha,\alpha'}( \mathcal{T}_\alpha ; (X_i)_{i \geq 0} )$ the closure of this  increasing union.

Notice that when $ \alpha'=2$ we have $p_{\alpha,2,d,3}=0$ and thus no branch point of degree larger or equal to $4$ is created in $ \tau_{n}$. In other words, $(\tau_{n})$ is a sequence a binary trees. In this case, the pruning procedure boils down to adding those points $X_{n+1}$ whose attachment to $\tau_{n}$ preserves the binary structure of the tree and thus $ \mathrm{Prun}_{\alpha,2}( \mathcal{T}_{\alpha} ; (X_{i})_{i\geq 0})$ can be seen as a deterministic function of $ \mathcal{T}_{\alpha}$ and of the leaves $X_{i}, i \geq 0$.

\begin{theorem}  \label{thm:pruning} The statement of Proposition \ref{thm:inside}
 holds with  \begin{eqnarray*}\mathfrak{T}_{\alpha,\alpha'} &=&  \mathrm{Prun}_{\alpha,\alpha'}\big( \mathcal{T}_\alpha ; (X_i)_{i \geq 0} \big).   \end{eqnarray*}
\end{theorem}

Moreover, since $ \alpha' \mapsto p_{\alpha,\alpha',d,d'}
$ is decreasing, we can couple the realizations of the pruned subtrees such that $ \alpha' \in (\alpha,2] \mapsto \mathrm{Prun}_{\alpha,\alpha'}\big( \mathcal{T}_\alpha ; (X_i)_{i \geq 0} \big)$ is decreasing inside $ \mathcal{T}_\alpha$. The pruning operation can also be viewed from a fragmentation point of view which sheds new light on the construction of $ \mathfrak{T}_{\alpha,\alpha'}$.

\bigskip

\noindent \textbf{Modification of the stable fragmentation.}  Fragmentation theory is a mathematical attempt to modelize the behavior of particles that undergo a splitting process, see \cite{Ber06} for an account on this field.  Any $\alpha$-stable tree $ \mathcal{T}_{\alpha}$ can be associated with a self-similar fragmentation process that describes the masses of the subtrees of $ \mathcal{T}_{\alpha}$ above a certain height. In other words, an $\alpha$-stable tree can be seen as the genealogical tree (in the sense of \cite{HM04})  of a pure-jump self-similar fragmentation process with index  $-\bar \alpha$,  whose dislocation measure $\nu_{\alpha}$ has been identified by Bertoin \cite{Ber02} in the Brownian case $\alpha=2$ and by  Miermont \cite{Mie03} in the general case (see Section \ref{sec:frag} for an expression of $\nu_{\alpha}$). All these dislocation measures are conservative, which means that  the mass is kept at each dislocation. 
In Section \ref{sec:frag} we introduce for $1<\alpha < \alpha' \leq 2$ a modification of the dislocation measure $\nu_{\alpha}$ by (roughly speaking) keeping randomly some of the fragments created by $\nu_{\alpha}$. When $\alpha'=2$, this  simply consists in  keeping only two fragments, chosen proportionally to their sizes. In the other cases the procedure is more complex and depends on the probabilities $p_{\alpha,\alpha',d,d'}$ introduced in (\ref{pdd}). Obviously, the resulting dislocation measure $\nu_{\alpha,\alpha'}$ is now \emph{dissipative} (some mass is lost at each dislocation). It is still possible to associate with it a random tree coding its genealogy  (see \cite{Steph}) and we have:

\begin{proposition} \label{thm:frag}The tree $ \mathfrak{T}_{\alpha,\alpha'}$ of Proposition \ref{thm:inside}
 can be seen as the genealogical tree of a pure-jump self-similar fragmentation of index $-\bar \alpha$  and dislocation measure $\nu_{\alpha,\alpha'}$.  \end{proposition}

Thus, in addition to being the genealogy of a canonical conservative fragmentation process, the stable tree of parameter $\alpha'$  is also, up to a random scaling, the genealogical tree of a dissipative fragmentation process with dislocation measure $\nu_{\alpha,\alpha'}$ and auto-similarity index $-\bar \alpha$ for  \textit{all} $\alpha \in(1,\alpha')$.
Besides,  as such a dissipative fragmentation tree, $ \mathfrak{T}_{\alpha,\alpha'}$ naturally carries a Malthusian measure whose  total  mass is  distributed, up to a deterministic scaling, as the random variable $M_{\alpha,\alpha'}$ appearing in the random scaling in Proposition \ref{thm:inside}. In Theorem \ref{thFrag} of Section \ref{sec:frag}, we reinforce Proposition \ref{thm:frag} by showing that this dissipative fragmentation tree endowed with its Malthusian measure is distributed as an $\alpha'$-stable tree with its natural uniform mass measure, up to a random scaling depending on $M_{\alpha,\alpha'}$.

\bigskip

\begin{figure}[!h]
 \begin{center}
 \includegraphics[height=3.3cm,width=13.5cm]{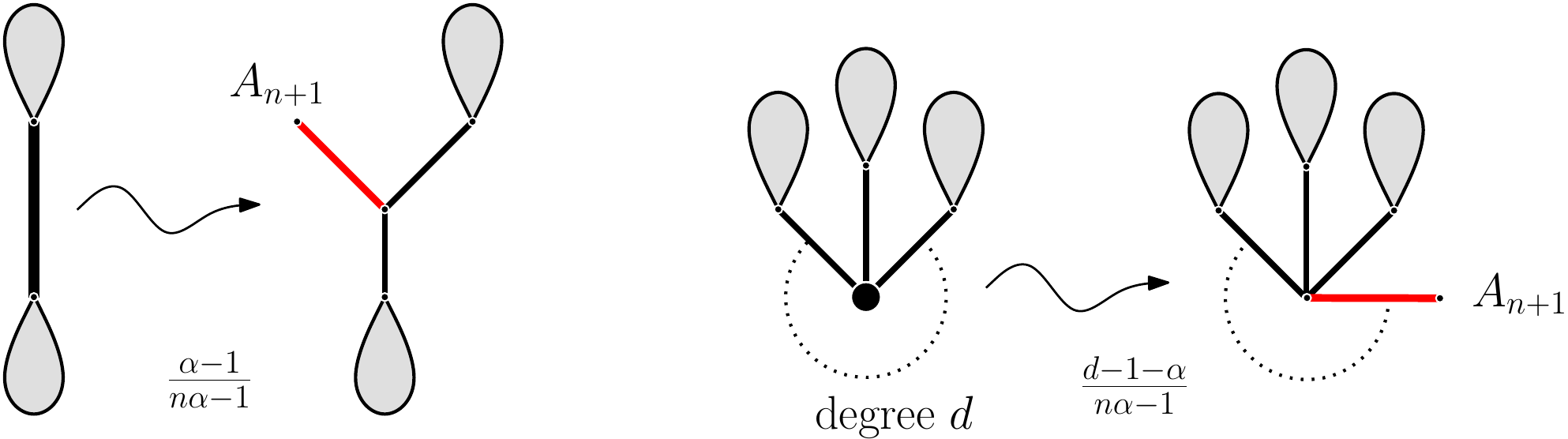}
 \caption{ \label{fig:marchal} illustration of the recursive rules to construct $ \mathbf{T}_{\alpha}(n+1)$ from $\mathbf{T}_{\alpha}(n)$.}
 \end{center}
 \end{figure}

\noindent \textbf{Strategy and organization of the paper.} The construction of the subtree $ \mathfrak{T}_{\alpha,\alpha'}$ presented in  Proposition \ref{thm:inside} relies on a discrete approach. In \cite{Mar08}, Marchal introduced a Markov chain  $( \mathbf{T}_{\alpha}(n), n \geq1)$ made of increasing labeled trees which, once re-normalized, converge towards the stable tree of parameter $\alpha$. Let us present quickly this construction: we start with $ \mathbf{T}_{\alpha}(1)$, the tree with a single edge, then inductively at step $n \geq 2$ we associate a weight $\alpha-1$ to every edges of the tree $ \mathbf{T}_{\alpha}(n)$ and a weight $d-1-\alpha$ to each vertex with degree $d \geq 3$. An edge or a vertex of $ \mathbf{T}_{\alpha}(n)$ is chosen at random accordingly to these weights and a new edge is attached either directly to the vertex or in the ``middle'' of the chosen edge, see Fig.\,\ref{fig:marchal} and Section \ref{SecMarchal} for more details. In the case $\alpha=2$ this construction is the famous algorithm for growing  binary trees due to  Rémy \cite{Rem85}.
We prove in Theorem \ref{thm:tool} that, once re-normalized, the trees $ \mathbf{T}_{\alpha}(n)$ converge towards a stable tree   \begin{eqnarray*}  \frac{ \mathbf{T}_{\alpha}(n)}{n^{\bar \alpha}} & \xrightarrow[n\to\infty]{\mathrm{a.s.}} & \alpha \mathcal{T}_{\alpha}, \end{eqnarray*} in the Gromov-Hausdorff sense. This result thus strengthens the  convergence in probability already obtained in \cite[Corollary 24]{HMPW08}. 

The key observation that triggered this work is  that  it is possible to identify within Marchal's construction of the $ \mathbf{T}_{\alpha}(n)$'s a ``sub-Markov chain'' of trees $ \mathbf{T}_{\alpha,\alpha'}(n) \subset \mathbf{T}_{\alpha}(n)$ whose growth mechanism is identical to Marchal's algorithm but with parameter $\alpha'>\alpha$, see Section \ref{SecDiscrete}. The scaling limit of the sequence $( \mathbf{T}_{\alpha,\alpha'}(n), n \geq 1)$ then  furnishes the tree $ \mathfrak{T}_{\alpha,\alpha'}$ of Proposition \ref{thm:inside}. \\

The paper is organized as follows. The first section contains the background on discrete and continuous trees as well as the presentation of Marchal's construction and its almost sure convergence towards the stable tree (Theorem \ref{thm:tool}). We then move in Section \ref{SecDiscrete}  to the observation that sub-constructions lie inside Marchal's algorithm and deduce Proposition \ref{thm:inside} and Theorem \ref{th:processus}. The third section is devoted to the pruning construction of $ \mathfrak{T}_{\alpha,\alpha'}$ (Theorem \ref{thm:pruning}). Finally, the last section explores the fragmentation approach and its consequences. \bigskip

\bigskip

\noindent \textbf{Acknowledgments:} this work is partially supported by  ANR-08-BLAN-0190 and ANR-08-BLAN-0220-01.

\clearpage

\section{Background on stable trees}
\label{SecStable}
In this section, we present the recursive construction of Marchal and prove (Theorem \ref{thm:tool}) that it converges \emph{almost surely}, after rescaling, towards a stable tree. Before embarking into that topic, we start by introducing some background on trees and the Gromov-Hausdorff and Gromov-Hausdorff-Prokhorov distances.

\subsection{Discrete and continuous trees}

\noindent \textbf{Discrete trees.} A discrete {tree} $\tau$ is a finite connected graph without cycle, considered  up to graph isomorphisms: it is not embedded in any space and its vertices are unlabeled.  If $x$ and $y$ are two vertices of a tree $\tau$, we denote by $ \llbracket x,y\rrbracket$ the discrete geodesic in $\tau$ between $x$ and $y$. If $y_{1}, y_{2},... , y_{k}$ are (distinct) vertices of $\tau$ we write 
$$ \mathrm{Span}(\tau; y_{1}, ... , y_{k}) =  \bigcup_{1 \leq i,j \leq k} \llbracket y_{i}, y_{j}\rrbracket, $$ for the discrete tree spanned by these vertices. The degree of a vertex is the number of edges adjacent to it, for instance a leaf is a vertex of degree one. A tree is binary if the degrees of its vertices are in $\{1,2,3\}$. 

A labeled discrete tree $\boldsymbol{\tau}=(\tau; x_0, x_1, ... , x_n)$ is a pair formed by a discrete  tree $\tau$ given with an exhaustive enumeration of its \emph{leaves}. If $\boldsymbol{\tau}$ is a labeled tree we call the {shape} of $ \boldsymbol{\tau}$ the tree $\tau$ obtained after forgetting the labeling of the tree. We will systematically use bold letters for labeled trees and standard ones for their associated shapes. 
\bigskip
 
\noindent \textbf{Continuous trees.}   We briefly recall here some facts on  $\mathbb R$-trees and refer to \cite{Eva08,LG06} for an overview on this topic.  
An $\mathbb{R}$-tree is a metric space $(\mathcal T,d)$ such that for every $
x,y\in \mathcal T$,
\begin{enumerate}
\item[-] there is an isometry $\varphi_{x,y}:[0,d(x,y)]\to \mathcal T$ such that $
\varphi_{x,y}(0)=x$ and $\varphi_{x,y}(d(x,y))=y$
\item[-] for every continuous, injective function $c:[0,1]\to \mathcal T$ with
  $c(0)=x,$ $c(1)=y$, one has $c([0,1])=\varphi_{x,y}([0,d(x,y)])$.
\end{enumerate}
We identify two  $\mathbb R$-trees when they are isometric, and still use the notation $(\mathcal T,d)$ to design an isometry class.
Note that a discrete tree may be seen as a $\mathbb R$-tree by ``replacing" its edges by segments. Unless specified, it is implicit in this paper that these segments are all of length 1.  We use much of the notation we introduced in the context of discrete trees when it is non ambiguously extended to $\mathbb R$-trees. In particular if $(\mathcal{T},d)$ is an $ \mathbb{R}$-tree and if $a,b \in \mathcal{T}$ we denote by $ \llbracket a,b\rrbracket$ for the geodesic line between $a$ and $b$ in $ \mathcal{T}$. Also $ \mathrm{Span}({ \mathcal{T}}; y_{1}, ... , y_{k}) =  \cup_{1 \leq i,j \leq k} \llbracket y_{i}, y_{j}\rrbracket$ still denotes the subtree spanned by $y_{1}, ... , y_{k} \in \mathcal{T}$.  The degree (or multiplicity) of a point $x \in \mathcal{T}$ is the number of connected components of $ \mathcal{T}\backslash \{x\}$. A leaf is a point of degree $1$, a branch point has degree larger than or equal to $3$, and an $\mathbb R$-tree is said to be binary if the degrees of its points are in $\{1,2,3\}$.

\bigskip

\noindent \textbf{Gromov-Hausdorff-Prokhorov topology.}  The reader interested by the Gromov-Hausdorff and Gromov-Hausdorff-Prokhorov topologies should consult \cite{BBI01,Eva08,miermont09} for details and proofs. Let $k \in \{0,1,2,...\}$. A compact metric space $(E,d)$ is $k$-pointed if it is given with $k$ points $x_{1}, ... , x_{k} \in E$ (when $k=0$ this is just a compact metric space). An isometry between two $k$-pointed compact metric spaces is an isometry between the spaces that maps the $k$ distinguished points to each other. The set of isometry classes of $k$-pointed compact metric spaces is endowed with the classical ($k$-pointed) Gromov-Hausdorff topology, that makes it Polish. 
This distance can be defined as
 \begin{eqnarray*} \mathrm{d_{GH}}( (E,d ; x_{1}, ..., x_{k}),(E',d'; x_{1}', ... , x_{k}')) &=& \inf \Big \{ \mathrm{d_{H}}(\phi(E),\phi(E')) \vee \max_{1\leq i \leq k} \mathrm{\delta}(\phi(x_{i}),\phi'(x_{i}'))\Big\},  \end{eqnarray*} where the infimum is taken over all choices of  metric spaces $(F, \delta)$ and isometric embeddings $\phi : E \to F$, $\phi' : E' \to F$, and  where $ \mathrm{d_{H}}$ denotes the Hausdorff distance in $F$.  The class of $\mathbb R$-trees forms the accumulation points of the class of rescaled discrete trees for the Gromov-Hausdorff topology and is thus closed. In the remainder of the paper, we use the notation $a \cdot E$ for the rescaled metric space $(E,ad)$ for any $a>0$.

  The Gromov-Hausdorff topology can be enriched in order to take into account measured spaces. A compact metric space $(E,d)$ endowed with a Borel probability measure $\mu$ is called a measured compact metric space. We can extend the Gromov-Hausdorff topology to isometry classes (defined in the obvious way) of measured compact metric spaces by putting 
  \begin{eqnarray*} \mathrm{d_{GHP}}((E,d,\mu),(E',d',\mu')) &=& \inf  \left\{\mathrm d_{\mathrm{H}}(\phi(E),\phi'(E'))\vee \mathrm d_{\mathrm{P}}(\phi_*\mu,\phi'_*\mu')\right\}\, , \end{eqnarray*}
where again $\phi,\phi'$ are isometries from $E,E'$ into a common space $(F,\delta)$, $\phi_*\mu,\phi'_*\mu'$ are the push-forward of
$\mu,\mu'$ by $\phi,\phi'$, and $\mathrm d_{\mathrm{P}}$ is the Prokhorov
distance between Borel probability measures on $F$: $$\mathrm d_{\mathrm{P}}(m,m')=\inf\{\varepsilon>0:m(C)\leq m'(C^\varepsilon)+\varepsilon\mbox{ for every }C\subset F\mbox{ closed}\}\, ,$$ where $C^\varepsilon=\{x\in F:\inf_{y\in C}\mathrm d(x,y)<\varepsilon\}$ is the $\varepsilon$-neighborhood of $C$. The function $\mathrm d_{\mathrm{GHP}}$ is a distance on the set of isometry classes of measured compact metric spaces that makes it Polish. 

Let us give another convenient way to express Gromov-Hausdorff distances. A correspondence between  two $k$-pointed compact metric spaces $(E ,d ; x_{1}, ... , x_{k})$ and $(E',d' ; x_{1}', ... , x_{k}')$ is a subset $\mathcal{R}$ of $E\times E'$  such that, for every $y_{1} \in E$, there exists at least one point $y_{2}\in E'$ such that $(y_{1},y_{2}) \in \mathcal{R}$ and conversely, for every $z_{2}\in E'$, there exists at least one point $z_{1}\in E$ such that $(z_{1},z_{2}) \in \mathcal{R}$. Also $(x_{1},x_{1}'), ... , (x_{k}, x_{k}') \in \mathcal{R}$. The distortion of the correspondence $\mathcal{R}$ is defined by    \begin{eqnarray*} \mathrm{dis}(\mathcal{R}) &=& \sup_{}\big\{|d(x_{1},y_{1})-d'(x_{2},y_{2})| : (x_{1},x_{2}),(y_{1},y_{2}) \in \mathcal{R} \big\}.  \end{eqnarray*} The Gromov-Hausdorff distance between $(E,d; x_{1}, ... , x_{k})$ and $(E',d'; x_{1}',... , x_{k}')$ can be expressed as half the infimum of the distortions of correspondences between $E$ and $E'$. A similar (but more involved) definition of $ \mathrm{d_{GHP}}$ via correspondences can be found in \cite{miermont09}.

\subsection{The stable trees}

We now give some additional background on stable trees and refer to \cite{DLG02,DLG05} for a complete account. Stable trees are random variables $ \mathcal{T}_{\alpha}$ for $\alpha \in (1,2]$ taking values in the set of isometry classes of compact  $\mathbb R$-trees that can be defined in various manners. For example, the tree $\mathcal T_{\alpha}$ can  be defined  by its finite-dimensional marginals (see \cite[Theorem 3.3.3]{DLG02}) or by its contour process. In particular, the $2$-stable tree can be identified with $\sqrt{2}$ times the Brownian tree $\mathcal T_{\mathrm{Br}}$ coded by a standard normalized excursion $ (\mathbf{e}_{t} : 0 \leq t \leq 1)$. Note also that the Brownian tree introduced by Aldous in \cite{Ald93} is equal to twice $ \mathcal{T}_{ \mathrm{Br}}$ so that we have in distribution $$ \mathcal{T}_{2} = \sqrt{2} \mathcal{T}_{ \mathrm{Br}} = \frac{ \mathcal{T}_{ \mathrm{Aldous}}}{\sqrt{2}}. $$
We mention that contrary to the Brownian tree where all the branch points are of degree $3$, in the stable case $1<\alpha<2$, they are of infinite multiplicity (\cite{DLG05}). 

Let us precise the definition via scaling limits of Galton-Watson trees and introduce the mass measure on $ \mathcal{T}_{\alpha}$. For $\alpha \in (1,2)$, consider a (unordered) Galton-Watson tree with offspring $\eta$ with mean 1 and such that $\eta(k)\sim Ck^{-1-\alpha}$ when $k\rightarrow \infty$. Let $T^{\mathrm{GW}}_n$ denote a version of this tree conditioned to have $n$ vertices, and equip it  with the uniform probability measure on its vertices, which is denoted by $\mu_n$.
Then by Duquesne \cite{Du03} (recall that $\bar \alpha=1-1/\alpha$),
 \begin{eqnarray*}
\left(\frac{T^{\mathrm{GW}}_n}{n^{\bar \alpha}},\mu_n \right) & \xrightarrow[n\to\infty]{( \mathrm{d})} & \left(\left(\frac{\alpha(\alpha-1)}{C\Gamma(2-\alpha)}\right)^{1/\alpha} \mathcal T_{\alpha}, \mu_{\alpha}\right), 
 \end{eqnarray*} for the Gromov-Hausdorff-Prokhorov topology, where $\mathcal T_{\alpha}$ is a stable tree of index $\alpha$ equipped with a probability measure $\mu_{\alpha}$ which can be  interpreted as the uniform mass measure on  $\mathcal T_{\alpha}$. A similar result holds for $\alpha=2$. It is well-known \cite{DLG02} that for all $\alpha \in (1,2]$, the measure $\mu_{\alpha}$ is actually fully supported by the set of leaves of $\mathcal T_{\alpha}$ and that it can be measurably constructed from $\mathcal T_{\alpha}$. It thus makes sense to speak about an i.i.d. sample $X_{n}, n \geq 0$ of leaves according to $ \mu_{\alpha}$ conditionally on $ \mathcal{T}_{\alpha}$. We shall repeatedly use the following easy fact: $\{X_{n}, n \geq 0\}$ is dense in $ \mathcal{T}_{\alpha}$ a.s.

\subsection{Marchal's recursive construction}
\label{SecMarchal}

Let $\alpha \in (1,2]$.  In \cite{Mar08}, Marchal proposed a recursive construction to build random finite trees that  converge in the scaling limit towards the $\alpha$-stable tree. The construction is in fact a Markov chain $( \mathbf{T}_{\alpha}(n))_{n\geq 1}$ with values in the set of labeled trees, such that $ \mathbf{T}_{\alpha}(n)$ has $n+1$ leaves and is defined as follows:  we start with the tree $ \mathbf{T}_{\alpha}(1)$ which is the only tree with one edge and two labeled vertices denoted by $A_{0}$ and $A_{1}$.  We then build recursively $ \mathbf{T}_{\alpha}(n+1)$ from $ \mathbf{T}_{\alpha}(n)$ by adding a new edge. More precisely, given $ \mathbf{T}_{\alpha}(n)$, assign a weight $\alpha-1$ to any edge of $ \mathbf{T}_{\alpha}(n)$ and a weight $d-1-\alpha$ to each vertex of degree $d \geq 3$ (all other vertices have zero weight). The total weight $W( \mathbf{T}_\alpha(n))$ of the tree $ \mathbf{T}_{\alpha}(n)$ is easily checked to be  \begin{eqnarray} \label{weight} W( \mathbf{T}_\alpha(n))&=& n\alpha-1,  \end{eqnarray} and is in particular independent of the shape of the tree. We then choose an edge or a vertex of $ \mathbf{T}_{\alpha}(n)$ proportionally to its weight:
\begin{itemize}
\item	if we picked an edge then we split it into two edges with a middle vertex on which we attach a new edge carrying the $n+1$th leaf denoted by $A_{n+1}$.
\item if a vertex has been selected then we attach a new edge carrying $A_{n+1}$ to it. \end{itemize}
 In \cite{Mar08}, Marchal exactly computed  the distribution of the  tree $ \mathbf{T}_{\alpha}(n)$. More precisely, he proved by induction on $n$ that if $ \boldsymbol{t}_{0}$ is a fixed labeled tree with $n+1$ leaves,  then  \begin{eqnarray*}\mathbb P ( \mathbf{T}_{\alpha}(n) = \boldsymbol{t}_{0} ) &=& \frac{ \prod_{v \in  \boldsymbol{t}_{0}} p_{ \mathrm{deg}(v)}}{ \prod_{i=1}^{n-1} (i\alpha-1)},  \end{eqnarray*}
   \begin{eqnarray*} \mathrm{with} \quad p_{1}=1, \quad p_{2}=0 \quad \mbox{and} \quad  \quad p_{k} = |(\alpha-1) (\alpha-2) ... (\alpha-k+2)| \quad \mbox{for }k \geq 3.  \end{eqnarray*} From this formula, it is easy to deduce (see \cite{Mar08}) that  $\mathbf{T}_{\alpha}(n)$ as the same distribution as  the labeled tree obtained from the finite $n+1$-dimensional marginal of the  standard $\alpha$-stable tree, see Theorem 3.3.3 in \cite{DLG02}. We this identification in hands, it is shown in \cite[Corollary 24]{HMPW08} that there exits a stable tree  $\mathcal T_{\alpha}$ built on the same probability space that supports $\mathbf T_{\alpha}(n)$, $n \geq 1$, and equipped with its uniform mass measure $\mu_{\alpha}$, such that
 \begin{eqnarray*}
\left(\frac{\mathrm T_{\alpha}(n)}{\alpha n^{\bar \alpha}}, \frac{\sum_{i=0}^n \delta_{A_i}}{n+1} \right)&\xrightarrow[n\to\infty]{(\mathbb P)} & \left( \mathcal T_{\alpha},  \mu_{\alpha}\right)
 \end{eqnarray*}
for the Gromov-Hausdorff-Prokhorov topology. The goal of the next section is to improve this convergence into an almost-sure one.

\subsection{The convergence theorem}

\begin{figure}[!h]
 \begin{center}
 \includegraphics[width=12cm]{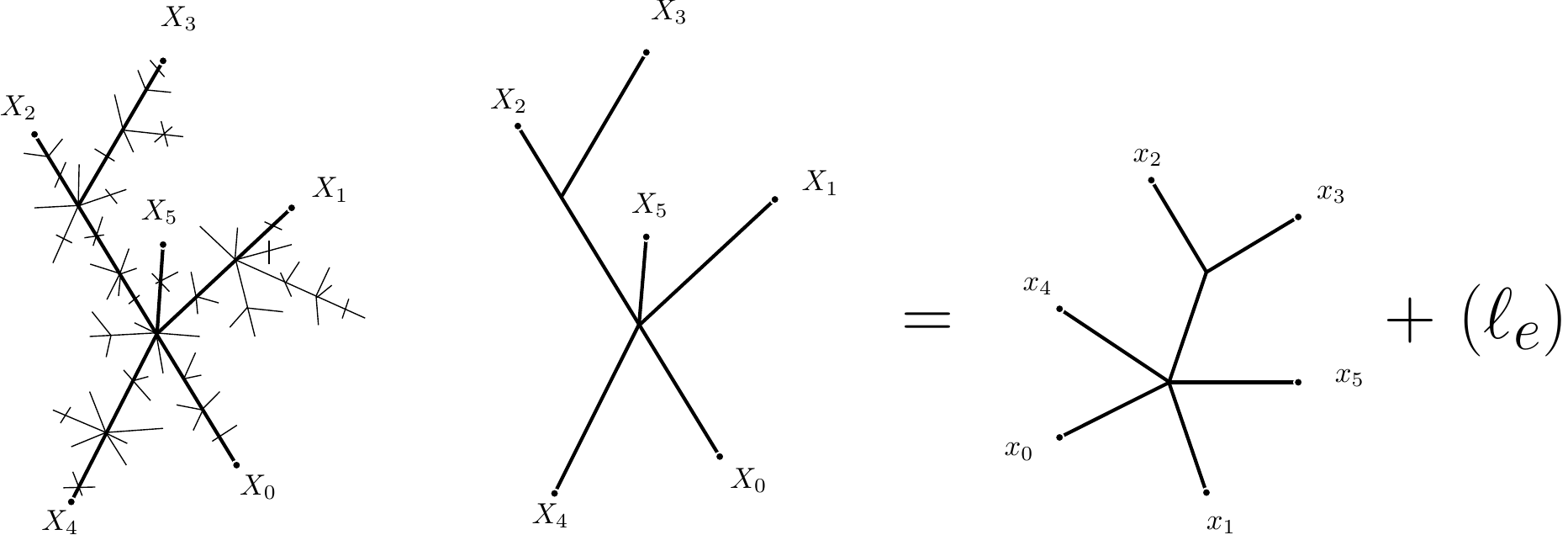}
 \caption{ \label{arbres} the trees $ \mathcal{T}_\alpha$, $ \mathrm{Span}( \mathcal{T}_\alpha ; X_0,..., X_5)$ and $\tau_{\alpha}(5)$ with its edge-lengths.}
 \end{center}
 \end{figure}

Let $\mathcal T_{\alpha}$ be an  $\alpha$-stable tree and conditionally on $ \mathcal{T}_{\alpha}$, let $(X_{i}, i \geq 0)$ be a sample of i.i.d. leaves distributed according to the mass distribution $\mu_{\alpha}$ on $ \mathcal{T}_{\alpha}$. We consider the finite $n+1$th-dimensional trees $\mathrm{Span}( \mathcal{T}_\alpha ; X_0,X_1, ... , X_n)$ for $n \geq 1$. These objects are by definition continuous trees but can equivalently be seen as discrete labeled trees $$ \boldsymbol{\tau}_{\alpha}(n)= (\tau_\alpha(n) ; x_0, x_1, ... , x_n)$$ carrying positive lengths $(\ell_{e})$ on their edges, see Fig.\,\ref{arbres}. 

\begin{theorem} \label{thm:tool} $(i)$ We have the joint identity in distribution
 \begin{eqnarray*}
\left( \mathbf{T}_{\alpha}(n),n\geq 1\right)&\overset{(\mathrm{d})}=&\left( \boldsymbol{\tau}_\alpha(n),n\geq 1\right).
 \end{eqnarray*}
$(ii)$ For every $k \geq 0$, we  have in the sense of $k$-pointed Gromov--Hausdorff distance
 \begin{eqnarray} 
  \left(\frac{\tau_{\alpha}(n)}{\alpha n^{ \bar{ \alpha}}} ; x_{0}, x_{1}, ... , x_{k}\right) & \xrightarrow[n\to\infty]{\mathrm{a.s.}}& (  \mathcal{T}_{\alpha} ; X_{0}, X_{1}, ... , X_{k}).   \label{eq:couplage}  \end{eqnarray}
  $(iii)$ If $\mu_{n} = \frac{1}{n+1} \sum_{i=0}^{n} \delta_{x_{i}}$ denotes the uniform mass measure on the leaves of $ \tau_{\alpha}(n)$ then we have the following almost sure convergence for the Gromov-Hausdorff-Prokhorov topology 
   \begin{eqnarray*}
\left(\frac{ \tau_{\alpha}(n)}{\alpha n^{{\bar{\alpha}}}}, \mu_{n}\right) & \xrightarrow[n\to\infty]{\mathrm{a.s.}}& ( \mathcal{T}_{\alpha}, \mu_{\alpha}).   \end{eqnarray*} 

 \end{theorem}
 
Consequently, $\mathrm T_{\alpha}(n)/ \alpha n^{\bar \alpha}$ endowed with the uniform mass measure on its leaves indeed converges almost surely for the Gromov-Hausdorff-Prokhorov topology towards an $\alpha$-stable tree endowed with its uniform mass measure.
 
 \proof To simplify the notation during the proof we fix $\alpha \in (1,2]$ once  for all and drop the indices $\alpha$ by writing $ \boldsymbol{\tau}_{n} = (\tau_{n} ; x_{0}, ... , x_{n})$ instead of $ \boldsymbol{\tau}_{\alpha}(n) = (\tau_{\alpha}(n) ; x_{0}, ... , x_{n})$. To get (i), we use the fact discussed at the end of the previous section, and proved by Marchal in \cite{Mar08}, that  $ \mathbf{T}_{\alpha}(n) = \boldsymbol{\tau}_n$ in distribution in terms of labeled trees, for all $n\geq 1$. Since  for $k \leq n$ the labeled trees $ \boldsymbol{\tau}_k$ and $ \mathbf{T}_{\alpha}(k)$ are deduced from $ \boldsymbol{\tau}_n$ resp.\,$ \mathbf{T}_{\alpha}(n)$ by the same deterministic procedure, we deduce that $( \mathbf{T}_{\alpha}(k))_{1 \leq k \leq n}$ has the same law as $ ( \boldsymbol{\tau}_k)_{1 \leq k \leq n}$. Since this holds for all $n \geq 1$ we indeed get $(i)$. 
 
Let us turn to the second point and remark first that for every $k \geq 1$ we have the almost sure convergence in the sense of pointed Gromov--Hausdorff metric 
 \begin{eqnarray} \Big( \mathrm{Span}( \mathcal{T}_\alpha ; (X_i)_{0\leq i \leq n}) ; X_{0}, X_{1}, ... , X_{k}\Big) & \xrightarrow[n\to\infty]{\mathrm{a.s.}}& \Big(  \mathcal{T}_{\alpha} ; X_{0}, X_{1}, ... , X_{k}\Big) \label{eq:easy}.  \end{eqnarray} This is easy to obtain and follows from the fact that the sequence of leaves $(X_i)_{i \geq 0}$ is almost surely dense in the compact tree $ \mathcal{T}_\alpha$. The convergence \eqref{eq:couplage} will be established by comparing the labeled tree $ \boldsymbol{\tau}_n$ to its continuous analogue $\mathrm{Span}( \mathcal{T}_\alpha ; (X_i)_{0\leq i \leq n})$ obtained by dressing it up with edge lengths. More precisely, we recall that the continuous tree $\mathrm{Span}( \mathcal{T}_\alpha; (X_i)_{0 \leq i \leq n})$ can be seen as the discrete labeled tree $ \boldsymbol{\tau}_n$ where each edge $e$ has been replaced by a Euclidean segment of length $\ell^{(n)}_e$, see Fig. \ref{arbres} below. 
 
For notation convenience we  write \begin{quote} $ \mathrm{d}^{(n)}_{ \mathrm{gr}}$ for the graph metric in $\tau_n$ (with edge lengths $1$) \\ $ \mathrm{d}_\ell^{(n)}$ for the metric on the vertices of $\tau_n$ associated to the edge-lengths $( \ell^{(n)}_e)$. \end{quote} In particular the two metrics live on the set of vertices of $\tau_n$. Also, for $0 \leq i,j \leq n$ the $ \mathrm{d}_{\ell}^{(n)}$ distance between $x_{i}$ and $x_{j}$ in $\tau_{n}$ is equal to the distance between $X_{i}$ and $X_{j}$ in $ \mathcal{T}_{\alpha}$.  Our goal is  to show that  \begin{eqnarray} \label{eq:dist} \sup_{a,b \in {\tau}_n} \left|\frac{ \mathrm{d}^{(n)}_{ \mathrm{gr}}(a,b)}{ n^{{\bar{\alpha}}}}- \alpha \mathrm{d}_\ell^{(n)}(a,b)\right| & \xrightarrow[n\to\infty]{\mathrm{a.s.}}&0.  \end{eqnarray}
Let us first show how to use the last display to complete the proof of the theorem. For that we use the definition of Gromov-Hausdorff distance via correspondences. Let $\delta$ denote the metric on $\mathcal T_{\alpha}$ and define a correspondence $ \mathcal{R}_{n}$ between the trees $(\tau_{n}, n^{- \bar{\alpha}} \mathrm{d}_{ \mathrm{gr}}^{(n)})$ and $ (\mathrm{Span}( \mathcal{T}_{\alpha} ; X_{0}, ... , X_{n}),\alpha \delta)$ by declaring that $a \in \tau_{n}$ and $x \in \mathrm{Span}( \mathcal{T}_{\alpha} ; X_{0}, ... , X_{n})$ are in correspondence if $x$ belongs to an ``edge'' of the continuous tree that is combinatorially adjacent to the associate vertex $a$ in ${\tau}_n$. Note that $x_{i}$ is in correspondence with $X_{i}$ for all $0 \leq i \leq n$. The distortion of this correspondence is bounded by   \begin{eqnarray*} \mathrm{dist}( \mathcal{R}_{n}) & \leq &\sup_{a,b \in {\tau}_n} \left|\frac{ \mathrm{d}^{(n)}_{ \mathrm{gr}}(a,b)}{n^{{\bar{\alpha}}}}-\alpha \mathrm{d}_\ell^{(n)}(a,b)\right| + 2 \alpha \sup \ell_{e}^{(n)}.  \end{eqnarray*} Then \eqref{eq:dist} shows that these distortions vanish as $n \to \infty$. We then use \eqref{eq:easy} to get  \eqref{eq:couplage}.   
  
It thus remains to prove \eqref{eq:dist}.  Fix $\varepsilon>0$. For any $n \geq 0$, we denote by $ \mathcal{L}_n = \sum \ell^{(n)}_e$ the total length of the tree $ \mathrm{Span}( \mathcal{T}_\alpha ; X_0, ... , X_n)$ and by $| {{\tau}} _n|$ the number of edges of the discrete  tree $ { \tau}_n$. In order to simplify notation we set   \begin{eqnarray*} \chi_{n} &=& \frac{| {\tau}_n|}{ \mathcal{L}_{n}}.  \end{eqnarray*}It follows from  \cite[Theorem 3.3.3]{DLG05} that conditionally on $ \boldsymbol{\tau}_n$ and  on $ \mathcal{L}_n$, the edge lengths $\ell^{(n)}_e$ are distributed uniformly on $ \mathbb{R}_+^{| {\tau}_n|}$ subject to the condition $ \sum \ell^{(n)}_e = \mathcal{L}_n$. In particular, conditionally on $ \boldsymbol{\tau}_n$ and on $ \mathcal{L}_n$, if $a$ and $b$ are two vertices of $ \boldsymbol{\tau}_n$ then  the distribution of the random variable $ \mathrm{d}^{(n)}_\ell(a,b)$ is given by first dividing uniformly the interval $[0, \mathcal{L}_n]$ into $| {\tau}_n|$ pieces and summing the first $ \mathrm{d}^{(n)}_{ \mathrm{gr}}(a,b)$-th ones.  We can then apply the following lemma to deduce
 \begin{eqnarray} \label{eq:eq1}   \mathbb {P} \left( \left.\left|\mathrm{d}^{(n)}_\ell(a,b) - \frac{\mathrm{d}^{(n)}_{ \mathrm{gr}}(a,b) }{\chi_{n}}\right| \geq  \varepsilon \mathrm{d}^{(n)}_{ \ell}(a,b)  \ \right| \boldsymbol{\tau}_n, \mathcal{L}_n \right) & \leq &  c_{ \varepsilon}^{-1}\exp\left({-c_\varepsilon  \mathrm{d}^{(n)}_{ \mathrm{gr}}(a,b)}\right),\\
 \label{eq:eq2} \mathbb  P \left(  \mathrm{d}_{\ell}^{(n)}(a,b) \geq \varepsilon \mid  \boldsymbol{\tau}_{n},   \mathcal{L}_{n} \right) &\leq& c_{ \varepsilon}^{-1}\exp\left(- c_{ \varepsilon} \frac{\chi_{n}}{\mathrm{d}^{(n)}_{ \mathrm{gr}}(a,b)} \right).\end{eqnarray}
 
 \begin{lemma} Let $\varepsilon>0$. There exists a constant $c_{\varepsilon}>0$ such that for all $n \geq 1$ and for all $L >0$, if $0=U_{0} < U_1 <  ... < U_{n-1} < U_{n}=L$ is a uniform splitting of the  interval $[0,L]$ then for any $0 \leq k \leq n$ we have  
 \begin{eqnarray*} \mathbb P\left( \left|U_k-\frac{kL}{n}\right| \geq  \varepsilon U_k \right) & \leq & c_{ \varepsilon}^{-1}\exp(-c_\varepsilon k), \\
 \mathbb P \left( U_{k} \geq \varepsilon \right ) & \leq &   c_{ \varepsilon}^{-1}\exp \left( -c_{\varepsilon} \frac{n}{L k} \right).\end{eqnarray*}
 \end{lemma}
 \proof By standard properties of the uniform splitting of the unit interval  we have $L^{-1}U_k = T_k/T_n$ in distribution where the $T_i$'s are the arrival times of a standard Poisson process. In particular for all $ \varepsilon>0$, there exists $c>0$ such that $$\mathbb P(|T_{i}-i| \geq \varepsilon i) \leq \exp(- ci)$$ for all $i \geq 0$ which easily implies the first inequality of the lemma. For the second inequality, one has, for all $k\geq 1$,  \begin{eqnarray*}
\mathbb  P( U_{k} \geq \varepsilon) &=& \mathbb P\left( T_{k} \geq \frac{ \varepsilon T_{n}}{L}\right) \\ & \leq & \mathbb P\left(T_{k} \geq \frac{ \varepsilon (1- \varepsilon) n}{L}\right) + \mathbb P\left( T_{n} <(1- \varepsilon)n\right) \\ &\leq&\mathbb  P\left( \mathcal{P}( \varepsilon(1- \varepsilon)nL^{-1}) \leq k\right) + \exp( -cn),  \end{eqnarray*} where $ \mathcal{P}(a)$ denotes a Poisson random variable of parameter $a >0$. Next, since the function $a \in (0,\infty) \mapsto \mathcal P(a)$ is stochastically increasing, and since, for $q \in \mathbb N$, $\mathcal P(qk)$ is distributed as the sum of $q$ independent Poisson random variables of parameter $k$, one has,
$$\mathbb P\left(  \mathcal{P}(  \varepsilon(1- \varepsilon) nL^{-1}) \leq k\right) \leq \mathbb P\left( \mathcal P\left(\left\lfloor \varepsilon (1- \varepsilon) n L^{-1}k^{-1} \right\rfloor k\right) \leq k \right)\leq \mathbb  P\left( \mathcal{P}(k) \leq k\right)^{\lfloor  \frac{\varepsilon (1- \varepsilon) n}{Lk}\rfloor}.$$ Putting the pieces together and using $\mathbb  P( \mathcal{P}(k) \leq k) \to 1/2$ as $k \to \infty$ we get the second inequality.\endproof
 \begin{lemma} \label{lem:astuce}We have the almost sure convergence  \begin{eqnarray*} \frac{\chi_{n}}{n^{{\bar{\alpha}}}} & \xrightarrow[n\to\infty]{\mathrm{a.s.}} & \alpha.  \end{eqnarray*}
 \end{lemma}
\proof It is possible to prove the last lemma by analysing  separately the processes $( | {\tau}_n|)_{n\geq 1}$ and $( \mathcal{L}_n)_{n\geq 1}$ and prove that $| {\tau}_n| \sim \alpha n$ and $\mathcal{L}_n \sim n^{1/\alpha}$ a.s. as $n \to \infty$. However we bypass these calculations by using a result of \cite[Section 5.4]{HMPW08}: we have 
 \begin{eqnarray} \label{eq:astuce}  \frac{\mathrm{d}^{(n)}_{ \mathrm{gr}}(x_0,x_1)}{n^{{\bar{\alpha}}}} & \xrightarrow[n\to\infty]{\mathrm{a.s.}} & \alpha \delta(X_{0},X_{1}).  \end{eqnarray}
 Recall that for any $i , j \geq 0$ and any $n \geq \mathrm{max}(i,j)$ we have $ \mathrm{d}^{(n)}_\ell(x_i,x_j) = \delta(X_{i},X_{j})$. Let $ \varepsilon >0$ and set  \begin{eqnarray*}A_{n} &=& \left\{ \left|\delta(X_{0},X_{1})- \frac{\mathrm{d}^{(n)}_{ \mathrm{gr}}(x_0,x_1)}{\chi_{n}} \right| \geq \varepsilon \delta(X_{0},X_{1}) \right\}.  \end{eqnarray*} Our goal is to show that the probability that $A_{n}$ is realized for infinitely many $n$'s is 0, thus proving that $\mathrm{d}^{(n)}_{ \mathrm{gr}}(x_0,x_1)/\chi_{n}$ almost surely converges towards $ \delta(X_{0},X_{1})$. Identifying the limit with \eqref{eq:astuce} then proves the statement of the lemma. To do this, let $\beta \in (0,\bar \alpha)$ and consider the events $B_{n} = \{ \mathrm{d}^{(n)}_{ \mathrm{gr}}(x_0,x_1) \geq n^{\beta}\}$. By \eqref{eq:astuce} (and the fact that $\delta(X_0,X_1)>0$ a.s.) we get that $\mathbb P(\cup_{k \geq 1} \cap_{n \geq k} B_n)=1$. Hence, \begin{eqnarray*} \mathbb{P} ( A_{n} \ i.o ) & = &   \mathbb{P}( A_{n} \cap B_{n} \ i.o), \end{eqnarray*} where $i.o.$ means ``is realized for infinitely many $n$'s''. Using \eqref{eq:eq1}, we get that $ \mathbb{P}(A_{n}\cap B_{n})$ is less than $ {c}_{\varepsilon}^{-1} \exp( - {c_{ \varepsilon}}n^{ \beta})$ hence by Borel-Cantelli $ \mathbb{P}( A_{n} \cap B_{n} \ i.o.) =0$. This completes the proof of the lemma. \endproof

Coming back to the proof of Theorem \ref{thm:tool} $(ii)$, for any $ n \geq 1$, conditionally on $ \boldsymbol{\tau}_{n}$ and on $ \mathcal{L}_{n}$, let us evaluate the probability that two vertices $a,b \in \tau_{n}$ are such that $$\left|\mathrm{d}^{(n)}_\ell(a,b) - \frac{\mathrm{d}^{(n)}_{ \mathrm{gr}}(a,b) }{\chi_{n}}\right| \geq  \varepsilon ( \mathsf{Diam} \vee 1) + \sqrt{\chi_{n}^{-1}}\qquad (*),$$  where $ \mathsf{Diam}$ is the diameter of $ \mathcal{T}_{\alpha}$, that is the maximal distance between any pair of points in $ \mathcal{T}_{\alpha}$. Remark that $\mathsf{Diam}$ bounds from above all the quantities $ \mathrm{d}_{\ell}^{(n)}(a,b)$.  We then split the cases according to the graph distance between $a$ and $b$: 
   \begin{itemize}
   \item if $ \mathrm{d}_{ \mathrm{gr}}^{(n)}(a,b) \geq \sqrt{\chi_{n}}$ then the conditional probability of $(*)$ is bounded by $ c_{ \varepsilon}^{-1}\exp(- c_{ \varepsilon} \sqrt{\chi_{n}})$ by \eqref{eq:eq1},
   \item if $ \mathrm{d}_{ \mathrm{gr}}^{(n)}(a,b) < \sqrt{\chi_{n}}$ then the conditional probability of $(*)$ is bounded by  $ c_{ \varepsilon}^{-1}\exp(- c_{ \varepsilon} \sqrt{\chi_{n}})$ by \eqref{eq:eq2}.\end{itemize} Noticing that there are deterministically less than $2n$ vertices in $ {\tau}_{n}$, we deduce that conditionally on $ \boldsymbol{\tau}_{n}$ and $ \mathcal{L}_{n}$, the probability that $(*)$ holds for two vertices $a,b \in {\tau}_{n}$ is less than $ c_{ \varepsilon}^{-1}4n^2 \exp(-c_{ \varepsilon} \sqrt{\chi_{n}})$. We then use Borel-Cantelli's lemma together with Lemma \ref{lem:astuce} to get that  \begin{eqnarray*} \sup_{a,b \in {\tau}_n} \left|\frac{ \mathrm{d}^{(n)}_{ \mathrm{gr}}(a,b)}{\chi_{n}}- \mathrm{d}_\ell^{(n)}(a,b)\right| & \xrightarrow[n\to\infty]{\mathrm{a.s.}}&0.  \end{eqnarray*}  Finally we use Lemma \ref{lem:astuce} again to replace $\chi_{n}$ by $\alpha n^{{\bar{\alpha}}}$ in the last display and get \eqref{eq:dist}. 

The last point of the theorem is obtained in the same spirit. Since conditionally on $ \mathcal{T}_{\alpha}$ the $X_{i}$'s are i.i.d. according to $\mu_{\alpha}$  it follows that  $ \mu_{\alpha}^{(n)}=\frac{1}{n+1} \sum_{i=0}^{n} \delta_{X_{i}}$ is the empirical measure associated to $\mu_{\alpha}$ and then  \begin{eqnarray*} \left( \mathcal{T}_{\alpha},  \mu_{\alpha}^{(n)} \right ) & \xrightarrow[n\to\infty]{ \mathrm{a.s.}} & ( \mathcal{T}_{\alpha}, \mu_{\alpha}),  \end{eqnarray*} in the Gromov-Hausdorff Prokhorov sense. Futhermore, since $ \mathrm{Span}( \mathcal{T}_{\alpha} ; X_{0}, ... , X_{n})$ is a subtree of $ \mathcal{T}_{\alpha}$ that converges in the Gromov-Hausdorff sense towards $ \mathcal{T}_{\alpha}$ the last display holds with $( \mathcal{T}_{\alpha},  \mu_{\alpha}^{(n)} )$ replaced by $( \mathrm{Span}( \mathcal{T}_{\alpha} ; X_{0}, ... , X_{n}),  \mu_{\alpha}^{(n)})$. Recall now that $\mathrm{Span}( \mathcal{T}_{\alpha} ; X_{0}, ... , X_{n})$ can be seen as the tree  $ \boldsymbol{\tau}_{n}$ where each edge has been replaced by a Euclidean segment of length $\ell_{e}^{(n)}$. This identification transports the measure $\mu_{\alpha}^{(n)}$ onto the discrete atomic measure $\mu_{n}$. Thus it suffices to prove that 
  \begin{eqnarray*} \mathrm{d_{GHP}} \left( \left( \tau_{n}, \frac{ \mathrm{d}_{ \mathrm{gr}}^{(n)}}{\alpha n^{ \bar{ \alpha}}}, \mu_{n}\right),\left( \tau_{n},  \mathrm{d}_{\ell}^{(n)}, \mu_{n}\right)\right) & \xrightarrow[n\to\infty]{}& 0.  \end{eqnarray*}
 But this again follows from \eqref{eq:dist}. Indeed, the obvious correspondence $ \mathcal{R} = \{ (x,x) : x \in \tau_n\}$ between $( \tau_{n}, \frac{ \mathrm{d}_{ \mathrm{gr}}^{(n)}}{\alpha n^{ \bar{ \alpha}}})$ and $( \tau_{n}, \mathrm{d}_{\ell}^{(n)})$ has a vanishing distortion as $n \to \infty$ by \eqref{eq:dist}. This means (see \cite{BBI01}) that for any $ \varepsilon>0$ and for all $n$ large enough  we can isometrically embed these two trees into a common metric space $(F_n, \delta_n)$ such that two elements in correspondence (that is the same vertex of $\tau_n$ in the two embeddings of the trees) are at $\delta_n$-distance less than $ \varepsilon$ from each other. The image measures $\nu_n$ and $\nu_n'$ of $\mu_n$ by these embeddings then obviously satisfy $ \nu_n( C) \leq \nu_n'(C^\varepsilon)$ and $ \nu'_n( C) \leq \nu_n(C^\varepsilon)$
  for all Borel $C \subset F_n$. This suffices to prove the claim and to finish the proof of the theorem.
\endproof


\section{The discrete approach}

\label{SecDiscrete}

This section is mainly devoted to the construction of the $\alpha'$-stable subtrees $\mathfrak T_{\alpha,\alpha'}$ of $\mathcal T_{\alpha}$, for $\alpha' > \alpha$, as stated in Proposition \ref{thm:inside}, and then to the proof  of Theorem \ref{th:processus}. In order to establish the key couplings of Marchal's constructions that will be used for the proofs, we introduce once and for all a sequence $U_{i}, i \geq 2 $ of i.i.d.\,random variables uniformly distributed on $(0,1)$ that are independent of all other variables introduced so far and later on.

\subsection{Coupling different Marchal's constructions}
\label{sec:RM}
\begin{figure}[!h]
 \begin{center}
 \includegraphics[width=13cm]{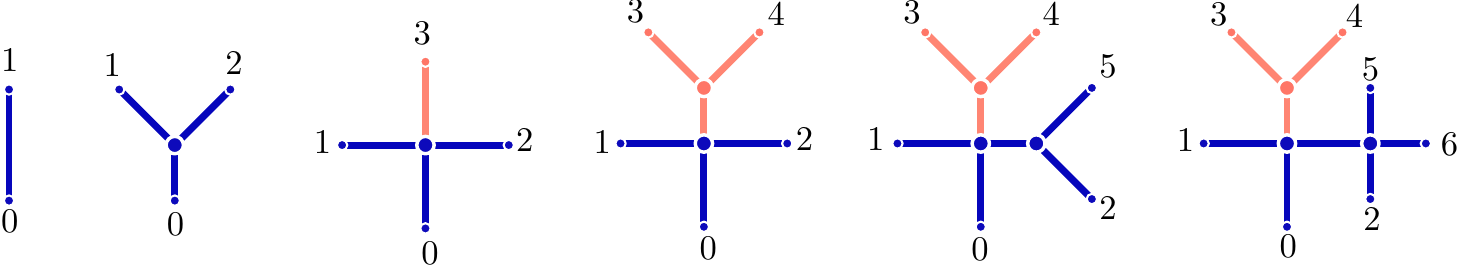}
 \caption{ \label{fig:insidebinary}illustration of the construction of the blue tree.}
 \end{center}
 \end{figure}

Fix $\alpha \in (1,2)$. We will couple Marchal's constructions for different values of $\alpha' \in (\alpha,2]$ in such a way that the $\alpha'$-constructions are  nested ``sub-constructions'' of the $\alpha$-one. To see this, fix $\alpha' \in (\alpha,2]$ and recall the definition of the probabilities $p_{\alpha,\alpha',d,d'}
$ in \eqref{pdd}, as well as the definition of Marchal's Markov chain $( \mathbf{T}_{\alpha}(n),n \geq 1)$ in Section \ref{SecMarchal}. The idea is to enrich the construction of this Markov chain with two colors, blue and red, such that the only edge and the two leaves of $ \mathbf{T}_{\alpha}(1)$ are blue and the colors evolve recursively according to the following rules, see Fig. \ref{fig:insidebinary}. At step $i \geq 2$:
\begin{enumerate}
\item[-] if an edge of color $ \mathsf{C} \in \{$red,blue$\}$ is selected, then we split this edge into two $ \mathsf{C}$ edges divided by a $ \mathsf{C}$ vertex and attach a $ \mathsf{C}$ edge-leaf to it;
\item[-] if a vertex is selected, then the new edge-leaf attached on it is blue if  $$U_{i} \leq p_{\alpha,\alpha',d_{i},d'_{i}}$$ and red otherwise, where $d_{i}$ is the degree of the vertex in the full tree $\mathbf T_{\alpha}(i-1)$ and $d'_{i}$ is its degree in the blue subtree of $\mathbf T_{\alpha}(i-1)$, with the convention that $d'_{i}=0$ if this vertex is red.
\end{enumerate} 

 We denote by $ \mathbf{T}_{\alpha,\alpha'}(n)$ the subtree of $ \mathbf{T}_{\alpha}(n)$ spanned by the blue edges and leaves. Remark that when $\alpha'=2$, the coloring rules are particularly simple (and deterministic):  each edge-leaf attached on a red edge or on a vertex (blue or red) is red, whereas an edge-leaf attached on a blue edge is blue. In this case the blue subtree is therefore binary. Remark also that the function $ \alpha' \in (\alpha,2] \mapsto p_{\alpha,\alpha',d,d'}$ is decreasing, hence  for every $n \geq 1$ the function $\alpha' \in (\alpha,2] \mapsto {T}_{\alpha,\alpha'}(n)$ is decreasing for the inclusion inside $ {T}_{\alpha}(n)$.

The labeling of $ \mathbf{T}_{\alpha,\alpha'}(n)$ is given by sorting the leaves of $ \mathbf{T}_{\alpha,\alpha'}(n)$ according to their labels in $ \mathbf{T}_{\alpha}(n)$. Also, we let $$L_{\alpha,\alpha'}(n)=\text{number of leaves of }  \mathbf{T}_{\alpha,\alpha'}(n) \text{ minus }1$$ and $L_{\alpha,\alpha'}^{-1}(n) = \inf\{ k \geq 1 : L_{\alpha,\alpha'}(k) =n\}$ its inverse. Here is the main observation:  
 \begin{lemma} \label{lemma:decoupling} The subchain $ \mathbf{T}_{\alpha,\alpha'}(.)$ evolves as a time-changed   Marchal's construction with parameter $\alpha'$. More precisely, we have the following identity in distribution   \begin{eqnarray*} \Big( \mathbf{T}_{\alpha,\alpha'}\big(L_{\alpha,\alpha'}^{-1}(n)\big)\Big)_{n \geq 1} & \overset{ \mathrm{(d)}}{=} & \big(\mathbf{T}_{\alpha'}(n)\big)_{n \geq 1}.  \end{eqnarray*}
 Futhermore, $( \mathbf{T}_{\alpha,\alpha'}(L_{\alpha,\alpha'}^{-1}(n)))_{n \geq 1}$ is independent of $(L_{\alpha,\alpha'}(n))_{n\geq 1}$.
 \end{lemma}
 
 \proof We examine the transition probabilities of the chain $( \mathbf{T}_{\alpha}(n), \mathbf{T}_{\alpha,\alpha'}(n))_{n \geq 1}$. Recall that  in Marchal's  construction, a weight $\alpha-1$ is assigned to each edge  and a weight $d-1-\alpha$ to each vertex of degree $d \geq 3$.  Without changing the dynamic, we multiply all these weights  by $(\alpha'-1)/(\alpha-1)$, so that each edge has now a weight $\alpha'-1$ and each vertex a weight $(d-1-\alpha)(\alpha'-1)/(\alpha-1)$. By \eqref{weight}, the total weight of the tree $ \mathbf{T}_{\alpha}(n)$ is  now
 $$ \overline{W}( \mathbf{T}_{\alpha}(n)) = \frac{\alpha'-1}{\alpha-1} (n \alpha-1).$$ 
 Hence, the new edge-leaf added to $ \mathbf{T}_{\alpha}(n)$ to get $ \mathbf{T}_{\alpha}(n+1)$ is blue if it has been
 \begin{enumerate}
 \item[-] grafted on an edge of $ \mathbf{T}_{\alpha,\alpha'}(n)$, which occurs with probability $(\alpha'-1)/ \overline{W}( \mathbf{T}_{\alpha}(n))$ for each edge of $ \mathbf{T}_{\alpha,\alpha'}(n)$
 \item[-] or grafted on a  vertex  of $ \mathbf{T}_{\alpha,\alpha'}(n)$ and then colored in blue, which occurs  with probability
 $$
\frac{(d-1-\alpha)(\alpha'-1)(\alpha-1)^{-1} \times p_{\alpha,\alpha',d,d'}}{\overline{W}( \mathbf{T}_{\alpha}(n)) }=\frac{d'-1-\alpha'}{\overline{W}( \mathbf{T}_{\alpha}(n))}
 $$
when $d'$ is the  degree of the vertex  in $ \mathbf{T}_{\alpha,\alpha'}(n)$ and $d$ its full degree  in $ \mathbf{T}_{\alpha}(n)$.  
\end{enumerate}
Observe that these probabilities are proportional to the  weights assigned  in Marchal's construction with parameter $\alpha'$. Consequently their sum, which is the total probability that the added edge-leaf is blue, is equal to 
$$
\frac{L_{\alpha,\alpha'}(n)\alpha'-1}{\overline{W}( \mathbf{T}_{\alpha}(n))},
$$
where we recall that $L_{\alpha,\alpha'}(n)+1$ is the number of leaves of $\mathbf{T}_{\alpha,\alpha'}(n).$
Moreover, conditionally  on the fact that the added edge-leaf is blue, it is grafted on an edge or vertex of  $ \mathbf{T}_{\alpha,\alpha'}(n)$  with the dynamic of Marchal's algorithm with parameter $\alpha'$.
In other words, the transition probabilities of the chain $( \mathbf{T}_{\alpha}(n), \mathbf{T}_{\alpha,\alpha'}(n))_{n \geq 1}$ can be described as follows. First,  $L_{\alpha,\alpha'}(n),n\geq 1$ is a Markov chain with transition probabilities
  \begin{eqnarray} \label{eq:chinese} \nonumber \mathbb P\left( L_{\alpha,\alpha'}(n+1)= L_{\alpha,\alpha'}(n)+1 \ | \ \mathcal F_n \right) &=& 1-\mathbb P\left(L_{\alpha,\alpha'}(n+1)=L_{\alpha,\alpha'}(n) \  | \ \mathcal F_n \right) \\[0.1cm]
  &=&\frac{L_{\alpha,\alpha'}(n)\alpha'-1}{\overline{W}( \mathbf{T}_{\alpha}(n))},
  \end{eqnarray} 
 where $\mathcal F_n$ denotes the sigma-field generated by $(\mathbf{T}_{\alpha}(k), \mathbf{T}_{\alpha,\alpha'}(k), k \leq n)$. 
 Second, conditionally on $L_{\alpha,\alpha'}(n+1)=L_{\alpha,\alpha'}(n)+1$, the transition probabilities on $ \mathbf{T}_{\alpha,\alpha'}(n)$ are those of  Marchal's construction with parameter $\alpha'$, whereas conditionally on $L_{\alpha,\alpha'}(n+1)= L_{\alpha,\alpha'}(n)$, the transition on the red component $\mathbf T_{\alpha}(n) \backslash  \mathbf{T}_{\alpha,\alpha'}(n)$ follows the appropriate conditional distribution. The statements of the lemma follow from these considerations. \endproof 

The chain $L_{ \alpha, \alpha'}$ can easily be studied using \eqref{eq:chinese}. Let us introduce its limit distribution (after re-normalization). A \emph{generalized Mittag-Leffler} random variable of parameters $(\beta,\theta)$ with $\beta \in (0,1)$ and $\theta > -\beta$ is a random variable denoted by $ \mathrm{ML}_{\beta,\theta}$ which is characterized by its positive moments \begin{equation}
\label{MomentsML}
\mathbb E\left[\mathrm{ML}^p_{\beta,\theta}\right]=\frac{\Gamma(\theta+1)\Gamma(\theta/\beta+p+1)}{\Gamma(\theta/\beta+1)\Gamma(\theta+p\beta+1)}, \quad p \geq 0.
\end{equation}
This random variable is actually a biased version of a power of a stable distribution of index $\beta$, see \cite[Section 3]{Pit06} for background. For later use, remark  that if $ 0 < a <b <1$ and if $\mathrm{ML}_{a,a}$, $ \mathrm{ML}_{b,b}$ and $ \mathrm{ML}_{a/b,a}$ are independent generalized Mittag-Leffler with parameters $(a,a)$, $(b,b)$ and $(a/b,a)$ respectively, then we have the following identity in distribution   \begin{eqnarray}  \label{eq:identity} \mathrm{ML}_{a,a} = \big(\mathrm{ML}_{a/b,a}\big)^b \cdot \mathrm{ML}_{b,b},  \end{eqnarray}this can be easily checked using moments \eqref{MomentsML}.

\begin{lemma}  \label{TnRRemi}There exists a generalized Mittag-Leffler r.v.  $ \mathrm{ML}_{\bar{\alpha}/\bar{\alpha}',\bar{\alpha}}$ such that
  \begin{eqnarray*} \frac{L_{\alpha,\alpha'}(n)}{n^{ \bar{\alpha}/\bar{\alpha}'}} & \xrightarrow[n\to\infty]{\mathrm{a.s.}} & \mathrm{ML}_{\bar{\alpha}/\bar{\alpha}',\bar{\alpha}}.  \end{eqnarray*}
  \end{lemma}
  \proof It is possible to analyse the chain $(L_{\alpha,\alpha'}(n))_{n\geq 1}$ ``by hand"  using the transition probabilities \eqref{eq:chinese}. However we bypass any calculation by using a connection with Pitman's Chinese restaurant process, see \cite[Section 3]{Pit06}. Indeed, it is straightforward to check from \eqref{eq:chinese} that the sequence $(L_{\alpha,\alpha'}(n+1)-1,n\geq 1)$ is distributed as the sequence of the number of tables in a $(\bar{\alpha}/\bar{\alpha}',\bar{\alpha})$-Chinese restaurant process, see \cite[Section 3.2.3]{Pit06}. The result then follows from \cite[Theorem 3.8]{Pit06}.
\endproof

\subsection{Proofs of Proposition \ref{thm:inside} and Theorem \ref{th:processus}} \label{sec:proofthem1}

Let $1<\alpha < \alpha' \leq 2$ and consider an $\alpha$-stable tree $\mathcal T_{\alpha}$ and its uniform mass measure $\mu_{\alpha}$. By Theorem \ref{thm:tool}, there exists a version of Marchal's Markov chain $(\mathbf T_{\alpha}(n), n \geq 1)$ such that 
 \begin{eqnarray} \label{cvmt}
\left(\frac{T_{\alpha}(n)}{\alpha n^{\bar{\alpha}}}, \frac{\sum_{i=0}^n \delta_{A_i}}{n+1} \right) & \xrightarrow[n\to\infty]{\mathrm{a.s.}} &  \left(\mathcal{T}_{\alpha}, \mu_{\alpha}\right),
 \end{eqnarray}
in the Gromov-Hausdorff-Prokhorov sense.
We will use the blue subchain $(\mathbf T_{\alpha,\alpha'}(n), n \geq 1)$  constructed from $(\mathbf T_{\alpha}(n), n \geq 1)$ in the previous section  to build the closed subtree $\mathfrak T_{\alpha,\alpha'}$ of Proposition  \ref{thm:inside}. 

\proof[Proof of Proposition \ref{thm:inside}]
Recall Lemma \ref{lemma:decoupling} and note from  Lemma \ref{TnRRemi} that $L_{\alpha,\alpha'}(n)\rightarrow \infty$ a.s. as $n \rightarrow \infty$. Then apply Theorem \ref{thm:tool}  to obtain the almost sure convergence  
\begin{eqnarray*} 
\frac{T_{\alpha,\alpha'}(n)}{{L_{\alpha,\alpha'}(n)}^{\bar{\alpha}'}} & \xrightarrow[n\to\infty]{\mathrm{a.s.}} &    \alpha' \mathcal{T}_{\alpha'},  \end{eqnarray*} in the Gromov-Hausdorff sense, where  $ \mathcal{T}_{\alpha'}$ denotes a stable tree of index $\alpha'$  (which is not independent of $\mathcal T_{\alpha}$!).
 Furthermore, still by Lemma \ref{lemma:decoupling}, the tree $  \mathcal{T}_{\alpha'}$ is independent of the sequence $(L_{\alpha,\alpha'}(n),n \geq 1)$, hence also of its limit after rescaling. 
Hence
\begin{eqnarray*}
\frac{T_{\alpha,\alpha'}(n)}{\alpha n^{{\bar{\alpha}}}} &=&  \frac{T_{\alpha,\alpha'}(n)}{\alpha {L_{\alpha,\alpha'}(n)}^{\bar{\alpha}'}}\times \left(\frac{L_{\alpha,\alpha'}(n)}{n^{\bar{\alpha}/\bar{\alpha}'}}\right)^{\bar{\alpha}'} \\
&\xrightarrow[n\to\infty]{\mathrm{a.s.}}& \frac{\alpha'}{\alpha} (\mathrm{ML}_{\bar{\alpha}/\bar{\alpha}',\bar{\alpha}})^{\bar{\alpha}'} \mathcal{T}_{\alpha'}=: \mathfrak{T}_{\alpha,\alpha'},
 \end{eqnarray*} where $\mathrm{ML}_{\bar{\alpha}/\bar{\alpha}',\bar{\alpha}}$ denotes the Mittag-Leffler random variable of Lemma \ref{TnRRemi}, which is independent of $ \mathcal{T}_{\alpha'}$. The variable $M_{\alpha, \alpha'}$ of Proposition \ref{thm:inside} is thus equal to $(\alpha'/\alpha)^{1/\bar \alpha '} \mathrm{ML}_{\bar{\alpha}/\bar{\alpha}',\bar{\alpha}}$ and the tree $\mathfrak T_{\alpha,\alpha'}$ can be realized as a closed subtree of $ \mathcal{T}_{\alpha}$. \endproof

In the constructions of $ \mathbf{T}_{\alpha}(n)$ and $ \mathbf{T}_{\alpha,\alpha'}(n)$, the two leaves $A_{0}$ and $A_{1}$ are kept in both trees and provide in the limit, by Theorem \ref{thm:tool}, two independent uniform points in $ \mathcal{T}_{\alpha}$ and in $ \mathfrak{T}_{\alpha,\alpha'}$. More precisely, for $\beta \in (1,2]$ we denote by $I_{\beta}$ the height between two uniform points in a $\beta$-stable tree. We also write $H_{n}$ the distance between $A_{0}$ and $A_{1}$ in $ \mathbf{T}_{\alpha}(n)$. Then we have by Theorem \ref{thm:tool} 
\begin{eqnarray*} \frac{H_{n}}{n^{ \bar{\alpha}}} & \xrightarrow[n\to\infty]{\mathrm{a.s.}} & \alpha I_{\alpha}.  \end{eqnarray*}
The distance $H_n$ is also the distance between $A_0$ and $A_1$ in $ \mathbf{T}_{\alpha,\alpha'}(n)$. Hence, by Lemmas \ref{lemma:decoupling} and \ref{TnRRemi}, combined with Theorem \ref{thm:tool},   \begin{eqnarray*}  \frac{H_{n}}{n^{{\bar{\alpha}}}}= \frac{H_{n}}{{L_{\alpha,\alpha'}(n)}^{\bar{\alpha}'}}\times \left(\frac{L_{\alpha,\alpha'}(n)}{n^{\bar{\alpha}/\bar{\alpha}'}}\right)^{\bar{\alpha}'} & \xrightarrow[n\to\infty]{\mathrm{a.s.}} & \alpha'I_{\alpha'} \left(\mathrm{ML}_{\bar{\alpha}/\bar{\alpha}',\bar{\alpha}}\right)^{\bar{\alpha}'},  \end{eqnarray*}
where, in the limit, the  Mittag-Leffler random variable and $I_{\alpha'}$ are independent. 
Thus if for $\alpha<\alpha' \in (1,2]$ we set 
 \begin{eqnarray*}  Q_{\alpha \to \alpha'} & \overset{ \mathrm{(d)}}{=}& \frac{\alpha'}{\alpha}\left(\mathrm{ML}_{\bar{\alpha}/\bar{\alpha}',\bar{\alpha}}\right)^{\bar{\alpha}'},  \end{eqnarray*} we can easily check that for $a<b<c \in (1,2]$,  we have  $ Q_{a \to c} = Q_{a \to b}\cdot Q_{b \to c}$ in distribution where the two last variables are independent. Hence, $Q$ can be interpreted as a transition kernel. Furthermore the random height in the stable trees form an invariant family for $Q$ in the sense that 
 \begin{eqnarray} I_{\alpha} &\overset{\mathrm{(d)}}=&   Q_{\alpha \to \alpha'}\cdot I_{\alpha'} \label{forward} \end{eqnarray} 
 with the two variables on the right-hand side independent.  This could have been proved directly using moments and using the fact (see e.g. \cite{Mie03}) that  $\alpha I_{\alpha}$ is distributed as $ \mathrm{ML}_{\bar{\alpha},\bar{\alpha}}$ for all $\alpha \in (1,2]$. Hence \eqref{forward} reduces to \eqref{eq:identity}.

We now aim at a backward analog of \eqref{forward}.  As in Theorem \ref{th:processus}, let $J_{\alpha}$ be distributed as the power of a Gamma distribution, namely $\alpha (\Gamma_{1+\bar \alpha})^{\bar \alpha}$  for $\alpha \in (1,2]$.
In particular the positive moments of $ J_{\alpha}$ are given by 
$$  \mathbb{E}\big[J_{\alpha}^p\big] =  \alpha ^p \frac{ \Gamma\left((p+1) \bar{\alpha}+1\right)}{ \Gamma\left( \bar{\alpha}+1\right)}, \qquad \mbox{for } p \geq 0. $$ 
These variables have been designed to satisfy the ``dual'' relations of \eqref{forward}, namely, for $1<\alpha<\alpha' \leq 2$, one can check using moments that the following identity in distribution holds \begin{eqnarray} \label{eq:id3} J_{\alpha}\cdot Q_{\alpha \to \alpha'} &=&  J_{\alpha'},  \end{eqnarray} when the two first random variables are independent. 
\bigskip

\noindent \textit{Proof of Theorem \ref{th:processus}.} Let $ \mathscr{T}_{\alpha}$ be a stable tree $ \mathcal{T}_{\alpha}$ that has been rescaled by an independent copy of $J_{\alpha}$, that is $ \mathscr{T}_{\alpha} = J_{\alpha} \cdot \mathcal{T}_{\alpha}$.  For $\alpha' \in (\alpha,2]$, using Proposition \ref{thm:inside},
 we can find inside $ \mathcal{T}_{\alpha}$ a subtree $ \mathfrak{T}_{\alpha,\alpha'}$ which is distributed as an $\alpha'$-stable tree rescaled by an independent factor distributed as $ Q_{\alpha \to \alpha'}$. Hence, after rescaling by $J_{\alpha}$ this provides a subtree of $ \mathscr{T}_{\alpha}$ which is distributed as 
$$ J_{\alpha} \cdot Q_{\alpha \to \alpha'} \mathcal{T}_{\alpha'}  \overset{ \mathrm{(d)}}{=} J_{\alpha'} \cdot \mathcal{T}_{\alpha'} \overset{ \mathrm{(d)}}{=} \mathscr{T}_{\alpha'},$$ since $ Q_{\alpha \to \alpha'}$ and $ J_{\alpha}$ are independent. Iterating, we obtain for all finite increasing sequences $1<\alpha_1<...<\alpha_n\leq 2$ a sequence of nested trees $ \mathscr{T}_{\alpha_n} \subset ... \subset  \mathscr{T}_{\alpha_1}$. Using Kolmogorov's extension theorem we can thus construct a process $ ( \mathscr{T}_{\alpha})_{1 < \alpha \leq 2}$ of nested rescaled stable trees. $\hfill \square$

\bigskip

Still using moments, note that 
$$
J_{\alpha}I_{\alpha} \overset{\mathrm{(d)}}=\Gamma_2,
$$
when $I_{\alpha}$ and  $J_{\alpha}$ are independent. This implies that the distance between two independent uniform points in $ \mathscr{T}_{\alpha}$ has a Gamma distribution of parameter 2 (that is with density $xe^{-x}$ on $(0,\infty)$) for all $\alpha \in (1,2]$. For the nested sequence of trees  $ ( \mathscr{T}_{\alpha})_{1 < \alpha \leq 2}$, it is then possible to couple the choice of the two   uniform independent points so that their distance is the \textit{same} in all trees $ \mathscr{T}_{\alpha}$.

\subsection{Zero measure}
\label{sec:zeromeasure}
\begin{proposition} 
\label{zeromass}
For all $\alpha' \in (\alpha,2]$, 
$$
\mu_{\alpha}\big( \mathfrak{T}_{\alpha,\alpha'}\big)=0 \quad \text{a.s}.
$$
\end{proposition}

\proof We will use that $(\mathcal T_{\alpha},\mu_{\alpha})$ is obtained from a Marchal's Markov chain $\mathbf T_{\alpha}$ via \eqref{cvmt}, jointly with the fact that $\mathfrak T_{\alpha,\alpha'}$ is the scaling limit of the subchain $\mathbf T_{\alpha,\alpha'}$. Since, by Lemma \ref{TnRRemi}, the number of leaves of $ \mathbf{T}_{\alpha,\alpha'}(n) $ is proportional to $n^{\bar{\alpha}/ \bar{\alpha}'} < \hspace{-1.5mm}< n$, it should be intuitively clear that $ \mathfrak{T}_{\alpha,\alpha'}$ has zero mass inside $ \mathcal{T}_{\alpha}$. We sketch here  a more rigorous argument. Let  $\gamma$ denote the expectation of the $\mu_{\alpha}$-mass of  $ \mathfrak{T}_{\alpha,\alpha'}$: $$ \gamma = \mathbb{E}\left[\mu_{\alpha} \big( \mathfrak{T}_{\alpha,\alpha'} \big)\right],$$ our goal being to show that $\gamma=0$. For that purpose we use the self-similarity of Marchal's construction. Indeed, if we stop this construction after one step,  we obtain the tree $\mathbf T_{\alpha}(2)$ made of a ``Y'' with three leaves $A_{0}, A_{1}$ and $A_{2}$ joined by a branch point denoted by $\mathcal B_{2}$. In the future evolution of the process, the three edges $e_{0}=\{\mathcal B_{2},A_{0}\}, e_{1}=\{ \mathcal B_{2},A_{1}\}$ and $e_{2}=\{\mathcal B_{2},A_{2}\}$ will give rise in the scaling limit to three copies of rescaled stable trees denotes by $ \tau_{0}, \tau_{1}$ and $\tau_{2}$. Similarly, the edges that will later be grafted on $\mathcal B_{2}$ will also give rise to a countable collection of continuous trees $ \tau_{3}, \tau_{4}, ...$ such that the tree $ \mathcal{T}_{\alpha}$ is made of the glueing of all the $ \tau_{i}$'s by one vertex. It is clear from Marchal's construction, and \eqref{cvmt}, that the measured trees $(\tau_i,\mu_{\alpha,i})$, $i\geq 0$, where $\mu_{\alpha,i}$ denotes the restriction of the measure $\mu_{\alpha}$ to $\tau_i$, satisfy $$\tau_i=\mu_{\alpha}(\tau_i)^{\bar \alpha}\mathcal T^{(i)}_{\alpha}, \quad \mu_{\alpha,i}=\mu_{\alpha}(\tau_i)\mu^{(i)}_{\alpha},$$ where the measured trees $(\mathcal T^{(i)}_{\alpha},\mu^{(i)}_{\alpha}), i \geq 0$ are i.i.d. copies of $(\mathcal T_{\alpha},\mu_{\alpha})$, that are moreover independent of the masses  $\mu_{\alpha}(\tau_i), i \geq 0$.

Let us now examine the evolution of the ``blue" subprocess $ \mathbf{T}_{\alpha,\alpha'}$ inside each of these trees. First, is it plain that the construction of $ \mathbf{T}_{\alpha,\alpha'}$ restricted to each offspring subtree of the three initial edges $e_{0},e_{1}$ and $e_{2}$ follows exactly the rules of a subchain of parameters $\alpha,\alpha'$ inside this subtree. This gives in the limit a subtree $t_{i} \subset \tau_{i}$, for $i\in \{0,1,2\}$ that satisfies   \begin{eqnarray} \label{frac} \gamma &=& \frac{\mathbb{E}[ \mu_{\alpha}( t_{i})]}{\mathbb{E}[ \mu_{\alpha}(\tau_{i})]}.  \end{eqnarray} For $i\geq 3$, the situation is almost the same, provided that the first edge grafted on $\mathcal B_{2}$ that will give rise to  $\tau_{i}$ is blue. If this is the case then the blue subprocess ``enters'' this part of the tree and provides in the limit a subtree $ t_{i} \subset \tau_{i}$ satisfying \eqref{frac}. If this is not the case, that is, the ancestral edge creating $\tau_{i}$ is red, then the blue subprocess does not enter this part of the tree and we set $t_{i} = \varnothing$. The subtree $ \mathfrak{T}_{\alpha,\alpha'}$ is thus made of the union of the $ t_{i},i\geq 0$ and we have   \begin{eqnarray} \label{exp} \gamma	 = \gamma  \mathbb{E}\Bigg[\sum_{i \geq 0} \mu_{\alpha}( \tau_{i}) \mathbf{1}_{\{{t}_{i} \ne \varnothing\}}\Bigg].  \end{eqnarray} To conclude, note that $\mathbb P(\exists i : t_i=\varnothing)=1$. Indeed, the probability that for all $\tau_i,i \geq 3$ the ancestral edges are blue is $$\prod_{d\geq 3} \frac{(d-1-\alpha')(\alpha-1)}{(d-1-\alpha)(\alpha'-1)}=0,$$ since $\alpha'>\alpha$. Together with the fact that $\sum_{i\geq 0} \mu_{\alpha}(\tau_i)=1$ a.s., this implies that the expectation in (\ref{exp}) is  strictly less than $1$. Hence $\gamma=0$.  \endproof

\endproof

\section{Pruning of discrete and continuous trees}

\label{sec:pruning} 
The goal of this section is to present a ``geometric" way of retrieving $ \mathbf{T}_{\alpha,\alpha'}(n)$ from the labeled tree $\mathbf{T}_\alpha(n)$. Although equivalent to the iterative construction of $ \mathbf{T}_{\alpha, \alpha'}(n)$ this new procedure, called ``pruning operation'' is well-suited to pass to the continuous limit  and yields Theorem \ref{thm:pruning}. 

\subsection{Pruning} \label{sec:pruning1} Let $1<\alpha<\alpha' \leq 2$. We recall and extend the the pruning procedure of stable trees described in the Introduction to the context of continuous or discrete trees. This procedure depends on the probabilities $p_{\alpha,\alpha',d,d'}$ and uses the sequence of uniform random variables $(U_{i})_{i \geq 2}$  introduced in Section 	\ref{SecMarchal}. Let $ ( {t} ; x_{0}, x_{1}, x_{2}, ... , x_{n})$ be a discrete or continuous tree given with an ordered subset of distinct leaves. We remove randomly some edges in this labeled tree so as  to obtain a subtree of $t$   (given with an ordered subset of its leaves) that we denote by $  \mathrm{Prun}_{\alpha,\alpha'}(t ; x_{0}, ... , x_{n})$ and which is constructed recursively as follows:
\begin{itemize}
\item the tree $\mathrm{Prun}_{\alpha,\alpha'}(t ; x_{0},x_{1})$ is $\llbracket x_{0},x_{1}	\rrbracket$, the tree spanned by $x_0$ and $x_1$ in $t$;
\item for $ i \geq 2$, consider two trees: $ \mathfrak{t}_{i-1}$, the tree spanned by $x_{0}, ... , x_{i-1}$ in $t$, and $\tau_{i-1}=\mathrm{Prun}_{\alpha,\alpha'}(t ; x_{0}, ... , x_{i-1})$. In $t$, the vertex $x_{i}$ is attached to $ \mathfrak{t}_{i-1}$ via the point $ \Delta_{i}$, that is $$\mathfrak{t}_{i-1} \cap \llbracket x_{i}, \Delta_{i}\rrbracket =\{\Delta_i\}.$$ Let then $d_{i}$ denote the degree of $\Delta_{i}$ in $ \mathfrak{t}_{i-1}$ and $d'_{i}$ its degree in $\tau_{i-1}$ with the usual convention that $d'_{i}=0$ if $\Delta_{i} \notin \tau_{i-1}$. See Fig. \ref{fig:binary}. We then set: $$
\begin{array}{ll} \mathrm{Prun}_{\alpha,\alpha'}(t ; x_{0}, ... , x_{i}) = \tau_{i-1} \cup \llbracket x_{i}, \Delta_{i}\rrbracket & \mbox{if } U_{i} \leq  p_{\alpha,\alpha',d_{i},d'_{i}},\\
 \mathrm{Prun}_{\alpha,\alpha'}(t ; x_{0}, ... , x_{i}) = \tau_{i-1} & \mbox{otherwise}.  \end{array}$$
\end{itemize} 
Note that $\mathrm{Prun}_{\alpha,\alpha'}(t ; x_{0}, ... , x_{n})$ is a function of $(t; x_{0}, ... , x_{n})$ and of the random variables $ (U_{i})_{i \geq 2}$. This tree can be given with an ordered sequence of leaves, which is the subsequence of elements of $\{x_{0}, x_{1}, ... , x_{n}\}$ that end-up in $ \mathrm{Prun}_{\alpha,\alpha'}(t ; x_{0}, ..., x_{n})$ and following our principles we denote this labeled tree by $ \mathbf{Prun}_{\alpha,\alpha'}(t ; x_0, ... , x_n)$. To stick with the color coding that we defined in Section \ref{sec:RM} we think of the pruned subtree as a blue subtree of $t$ whose complement in $\mathfrak{t}_n$ is red. In particular a leaf $x_i$ is blue if it belongs to the pruned subtree. All these constructions are coupled  \emph{via} the $U_{i}'$s in such a way that $$\mathrm{Prun}_{\alpha,\alpha'}(t ; x_{0}, ... , x_{n-1}) \subset \mathrm{Prun}_{\alpha,\alpha'}(t ; x_{0}, ... , x_{n}), \quad \forall n \geq 1,$$ and for every $\alpha$ the map $ \alpha'\in (\alpha,2] \mapsto \mathrm{Prun}_{\alpha,\alpha'}(t ; x_{0}, ... , x_{n}) $ is decreasing for the inclusion in $t$. Thus, when $ t$ is a continuous tree and $x_{0}, ... , x_{n}, ...$  is an infinite sequence of (distinct) leaves in $t$, we can define \begin{eqnarray*} \mathrm{Prun}_{\alpha,\alpha'}( t ; (x_{i})_{i \geq 0}) &=& \overline{\bigcup_{ n \geq 1} \mathrm{Prun}_{\alpha,\alpha'}( t ; x_{0}, ... , x_{n})}.  \end{eqnarray*}
Note also that when $ \alpha'=2$, since $ p_{\alpha,2,d,3}=0$, the tree $ \mathrm{Prun}_{\alpha,\alpha'}(t; x_{0}, ... , x_{n})$ is binary for all $n\geq 1$. In this case, the pruning procedure is deterministic (it does not depend on the $U_i$'s): the leaves $x_{i}$ are just added one by one provided that the tree remains binary. Notice also that the pruning constructions are coupled with the chain $\mathbf{T}_{\alpha,\alpha'}$ via the $U_i$'s.\bigskip

Now, recall that $\mathbf{T}_\alpha(n)$ is a discrete tree with $n+1$ ordered leaves $A_0,...,A_n$. 
The pruning has been defined such that the blue subtree $ \mathbf{T}_{\alpha,\alpha'}(n)$ can also be seen as $ \mathbf{Prun}_{\alpha,\alpha'}(\mathbf{T}_\alpha(n))$:

\begin{proposition} \label{prop:discretepruning} For all $ 1< \alpha < \alpha' \leq 2$ and for every $n \geq 1$,  we have    \begin{eqnarray*}\mathbf{Prun}_{\alpha,\alpha'}\big( \mathbf{T}_\alpha(n)\big) & {=}&   \mathbf{T}_{\alpha,\alpha'}(n).  \end{eqnarray*}
\end{proposition}

\proof This is easy to prove by induction on $n$ and is safely left to the reader. 
\endproof

\subsection{Proof of Theorem \ref{thm:pruning}}
\label{secproofpruning}

We start  with a continuity property of the application $ \mathrm{Prun}(.)$ with respect to the leaves. It applies both in the discrete and continuous settings and will be useful to prove Theorem \ref{thm:pruning}. 
\begin{proposition} \label{prop:lipschitz} Let $(t ; x_{0}, ... , x_{n})$ be a discrete or continuous tree given with a subset of $n+1$ leaves. Then for any $k \in \{0,1,..., n\}$ we have  \begin{eqnarray*}    \mathscr{D} \big( \mathrm{Prun}_{\alpha,\alpha'}(t ; x_0, ..., x_n) , \mathrm{Prun}_{\alpha,\alpha'}( t ; x_{0},... , x_{k})\big) &\leq & \inf_{\varepsilon>0} \left\{ x_{0}, ..., x_{k} \mbox{ is an } \varepsilon \mbox{-net in } (t, \mathscr{D})\right\},\end{eqnarray*}  where $ \mathscr{D}$ denotes either the graph distance if $t$ is discrete or its metric if $t$ is a continuous tree.
\end{proposition}
\proof Let $ 0 \leq k \leq n$. By the monotonicity of the pruning operation, the pruned tree $ \mathrm{Prun}_{\alpha,\alpha'}(t ; x_0, ..., x_n)$ can be obtained from $ \mathrm{Prun}_{\alpha,\alpha'}(t; x_{0}, ... , x_{k})$ by grafting on each of its vertices some trees.  The maximal height of these grafted trees is thus equal to the (Hausdorff) distance between $ \mathrm{Prun}_{\alpha,\alpha'}(t; x_{0}, ... , x_{k})$ and $ \mathrm{Prun}_{\alpha,\alpha'}(t; x_0, ... , x_n)$. We just have to remark that none of the points $x_{0}, ... , x_{k}$ belong to the grafted trees that is 
 \begin{eqnarray*} \{x_{0}, ... , x_{k}\} \cap \big( \mathrm{Prun}_{\alpha,\alpha'}(t; x_0, ..., x_n)\backslash  \mathrm{Prun}_{\alpha,\alpha'}( t; x_{0},... , x_{k}) \big) &=& \varnothing.  \end{eqnarray*} 
The statement of the proposition follows. \endproof

\proof[Proof of Theorem \ref{thm:pruning}] Starting with an $\alpha$-stable tree $\mathcal{T}_{\alpha}$  and a sample $(X_{i}, i \geq 0)$ of i.i.d. random {leaves} with distribution $\mu_{\alpha}$ (given $\mathcal T_{\alpha}$), we consider the chain  $ \mathbf{T}_{\alpha}$ such that (\ref{cvmt}) holds, as well as its blue subchain $\mathbf T_{\alpha,\alpha'}$.
We have already established (see the proof of Proposition \ref{thm:inside}) that once re-normalized by $ \alpha n^{\bar{\alpha}}$ the blue tree $ {T}_{\alpha,\alpha'}(n)$ converges almost surely towards $ \mathfrak{T}_{\alpha,\alpha'} \subset \mathcal{T}_{\alpha}$.  By Proposition \ref{prop:discretepruning} we just have to show that $ \mathrm{Prun}_{\alpha,\alpha'}( \mathbf{T}_{\alpha}(n))$ converges after renormalization towards $ \mathrm{Prun}_{\alpha,\alpha'}( \mathcal{T}_{\alpha} ; (X_{i})_{i \geq 0}).$ For this, we will use that for every $k \geq 0$ (see Theorem \ref{thm:tool} $(ii)$)   \begin{eqnarray} 
\label{cvkplus1}
\left( \frac{T_{\alpha}(n)}{\alpha n^{\bar{ \alpha}}}; A_{0}, A_{1}, ... , A_{k}\right) & \xrightarrow[n\to\infty]{\mathrm{a.s.}} & ( \mathcal{T}_{\alpha} ; X_{0}, ... , X_{k}),\end{eqnarray}
in the $k+1$-pointed Gromov-Hausdorff topology.  Let $ \varepsilon>0$.  Almost surely, since the $X_{i}$'s form a dense subset of $\mathcal{T}_{\alpha}$, which is compact, we can  find a (random) $0 \leq k < \infty$ such that $\{X_{0}, ... , X_{k}\}$ is an $ \varepsilon$-net in $ \mathcal{T}_{\alpha}$.  In particular, we eventually have that $\{A_{0}, ... , A_{k}\}$ is a $\lfloor 2 \alpha \varepsilon n^{{\bar{\alpha}}} \rfloor $-net in $ T_{\alpha}(n)$. Applying Proposition \ref{prop:lipschitz}  twice we get that for $n$ large enough, \begin{eqnarray*} \delta \big(\mathrm{Prun}_{\alpha,\alpha'}( \mathcal{T}_{\alpha} ; (X_{i})_{i \geq 0}),\mathrm{Prun}_{\alpha,\alpha'}( \mathcal{T}_{\alpha} ; X_{0},... , X_{k})  \big) & \leq & \varepsilon,\\
 \frac{1}{\alpha n^{{\bar{\alpha}}}}\mathrm{d_{gr}}\big(\mathrm{Prun}_{\alpha,\alpha'}( \mathbf{T}_{\alpha}(n)), \mathrm{Prun}_{\alpha,\alpha'}( {T}_{\alpha}(n) ; A_0, ... , A_k) \big) & \leq & 2 \varepsilon,\end{eqnarray*} 
 where $\delta$ is the metric in $\mathcal T_{\alpha}$.
Hence the proof will be completed when we will have shown  that for every fixed $k_{0} \geq 0$,   \begin{eqnarray} \label{eq:goal1}\frac{1}{ \alpha n^{{\bar{\alpha}}}} \mathrm{Prun}_{\alpha,\alpha'}({T}_{\alpha}(n) ; A_0, ... , A_{k_{0}})  & \xrightarrow[n\to\infty]{\mathrm{a.s.}}& \mathrm{Prun}_{\alpha,\alpha'}( \mathcal{T}_{\alpha} ; X_0, ..., X_{k_{0}}),  \end{eqnarray}  in the Gromov--Hausdorff sense as $n \to \infty$. 
To see this, we will obviously use (\ref{cvkplus1}). 
 Let us however emphasize here that the convergence \eqref{cvkplus1} in the $k+1$-pointed Gromov-Hausdorff sense \emph{does not imply} the convergence of the pruned subtrees in general. Counter-examples can easily be set up when some of the marked points in the limiting tree are not leaves. However, here, the $X_i$'s are distinguished leaves of $\mathcal T_{\alpha}$, therefore the discrete labeled trees obtained by removing vertices of degree 2 (and glueing the adjacent edges) in $\mathrm{Span}(T_{\alpha}(n);A_0,...,A_{k_0}), n \geq k_0$ and forgetting lengths on the edges of  $\mathrm{Span}(\mathcal T_{\alpha};X_0,...,X_{k_0})$ are all \textit{identical}. From this and the definition of the pruning procedure (with, for reminder,  the same sequence of uniform random variables $U_i,i \geq 2$ in all cases), we get that the discrete labeled trees obtained by removing vertices of degree 2 in $\mathrm{Prun}_{\alpha,\alpha'}(T_{\alpha}(n);A_0,...,A_{k_0})$ for $n$ large enough  and forgetting lengths on the edges of $\mathrm{Prun}_{\alpha,\alpha'}(\mathcal T_{\alpha};X_0,...,X_{k_0})$ are all identical. Re-incorporating distances, the conclusion follows from the convergence  of $\mathrm{Span}(T_{\alpha}(n);A_0,...,A_{k_0})$, after normalization,  towards  $\mathrm{Span}(\mathcal T_{\alpha};X_0,...,X_{k_0})$, which is a consequence of (\ref{cvkplus1}), see e.g. Lemma 14 in \cite{HM12}.
\endproof  

With the notation of the proof of Proposition \ref{thm:inside} we have thus proved that $$ \mathrm{Prun}_{\alpha,\alpha'}( \mathcal{T}_\alpha ; (X_i)_{i\geq0})=\mathfrak{T}_{\alpha,\alpha'}= M_{\alpha,\alpha'}^{\bar \alpha'} \cdot \mathcal{T}_{\alpha'}.$$ 
With a slight abuse of notation, we also denote by $\mu_{\alpha'}$ the uniform mass measure on the tree $ \mathrm{Prun}_{\alpha,\alpha'}( \mathcal{T}_\alpha ; (X_i)_{i\geq0})$. 
Let then $\mathsf{B}_i,  i \geq 0$ be the indices of the leaves that are kept in the pruning procedure applied to $ ( \mathcal{T}_{\alpha} ; (X_i)_{i \geq 0})$. By the above proof, these indices also correspond to the indices of the leaves $A_i, i \geq 0$ that end-up being blue in the construction of $ \mathbf{T}_{\alpha,\alpha'}$ from $ \mathbf{T}_{\alpha}$. We then claim that conditionally on  $\mathrm{Prun}_{\alpha,\alpha'}( \mathcal{T}_\alpha ; (X_i)_{i\geq0})$,  \begin{eqnarray*}   \mbox{ the leaves } X_{ \mathsf{B}_i}, i \geq 0 \mbox{ are i.i.d. according to }\mu_{\alpha'}.  \end{eqnarray*} In words, the blue leaves of the pruning procedure form an i.i.d.\,sequence of leaves distributed according to the uniform mass measure on the remaining pruned subtree. Although this is a purely continuous setup statement we give a proof relying on a discrete argument. Indeed, it follows from the preceding proof that for every $k \geq 0$ 
 \begin{eqnarray*} \bigg(\frac{T_{\alpha,\alpha'}(n)}{\alpha n^{\bar{\alpha}}} ; A_{ \mathsf{B}_0}, ... , A_{ \mathsf{B}_k} \bigg) & \xrightarrow[n\to\infty]{ \mathrm{a.s.}} &  \big(\mathrm{Prun}_{\alpha,\alpha'}( \mathcal{T}_\alpha ; (X_i)_{i\geq0}) ; X_{ \mathsf{B}_0}, ... , X_{ \mathsf{B}_k} \big),  \end{eqnarray*} in the $k+1$-pointed Gromov-Hausdorff sense. However, since $ \mathbf{T}_{\alpha,\alpha'}$ has the same distribution has a Marchal chain of parameter $\alpha'$ time-changed by an independent increasing process, we deduce from Theorem \ref{thm:tool} $(ii)$ that $ X_{ \mathsf{B}_0}, ... , X_{ \mathsf{B}_k}$ are i.i.d. according to $ \mu_{\alpha'}$ conditionally on $ M_{\alpha,\alpha'}^{\bar \alpha'}\cdot \mathcal{T}_{\alpha'}$ as demanded.

\section{A fragmentation point of view}
In this section, we exploit the pruning operation of the last section  in order to give another point of view on the construction of 
$ \mathfrak{T_{\alpha,\alpha'}}$ which is based on the fragmentation properties of the stable tree. Indeed, it is by now standard that the stable tree can be constructed from a certain self-similar fragmentation process whose conservative dislocation measure is explicitly known. We will show that we can trim this fragmentation process by throwing certain fragments away and such that the genealogy of the resulting \emph{dissipative} fragmentation process is coded by the tree $ \mathfrak{T}_{\alpha,\alpha'}$. This, finally, gives another interpretation of the random scaling appearing in Proposition \ref{thm:inside} as the limit of the Malthusian martingale associated to the dissipative fragmentation. 

\label{sec:frag}

\subsection{Fragmentation of the stable tree}
\label{sec:backgroundfrag}

\noindent \textbf{General self-similar fragmentation processes.} This paragraph is deliberately brief and we refer to Bertoin \cite{Ber01,Ber02,Ber06} for a  detailed and rigorous introduction to the theory of self-similar fragmentations. Let $a \in \mathbb R$ and $\nu$ be a sigma-finite measure on the set of decreasing sequences 
$$
\mathcal S=\Bigg\{\mathbf s=(s_1,s_2,...):s_1 \geq s_2 \geq ... \geq 0 \text{ and }\sum_{i\geq 1}s_i \leq 1  \Bigg\}
$$
(endowed with the topology of pointwise convergence),
such that $$\nu(1,0,...)=0 \quad \text{and} \quad \int_{\mathcal S}(1-s_1)\nu(\mathrm d \mathbf s)<\infty .$$ Such measure is called a \textit{conservative dislocation measure} when $\nu(\sum_i s_i<1)=0$ and a \textit{dissipative dislocation measure} otherwise.

A pure-jump self-similar fragmentation process $F$  with index of self-similarity $a$ and dislocation measure $\nu$ is a  $\mathcal S$-valued càdlàg Markov process modelizing the evolution of a system of splitting masses. The sequence $F(t)$ represents  the masses present at time $t$, ranked in the decreasing order.  When $\nu$ is finite the dynamic of the process is easy to describe: 
\begin{enumerate}
\item[(i)] different masses present at a given time evolve independently
\item[(ii)] each mass $m$ splits after a random time with an exponential distribution of parameter $m^a\nu(\mathcal S)$ to give sub-masses  $mS_1,mS_2,...$ where $(S_1,S_2,...) \in \mathcal S$ is distributed according to $\nu/\nu(\mathcal S)$, independently of the splitting time. 
\end{enumerate}
When $\nu$ is infinite, the splitting times are dense in $\mathbb R_+$ and the description of the process is more tricky. Informally, the item (i) is still valid, whereas the second item is replaced by the more general fact  that a mass $m$ splits in masses $m \mathbf s$, $\mathbf s \in \mathcal S$ at rate $m^a \nu(\mathrm d \mathbf s)$. 

In all cases, this can be made rigorous by using  Poisson point processes (see \cite{Berest02,Ber01} for details). 
We start with the description of \textit{homogeneous} fragmentations, the fragmentations with an index of self-similarity equal to 0.
Let $(S(t),k(t))_{t \geq 0}$ be a Poisson point process with intensity measure $\nu \otimes \#$, where $\#$ denotes the counting measure on $\mathbb N$. Then it is possible to build a pure-jump càdlàg $\mathcal S$-valued process $H$ such that $H(0)=(1,0,...)$ and $H$ jumps only when an atom $(S(t),k(t))$ of the Poisson point process occurs:  $H(t)$ is then obtained from $H(t-)$ by keeping all elements of the sequence $H(t-)$ but $H_{k(t)}(t-)$ which is  replaced by the sequence of masses $H_{k(t)}(t-)S_i(t), i \geq 1$, that is
$$
H(t)=\ \downarrow \big\{H_{j}(t-), H_{k(t)}(t-)S_i(t),i,j \geq 1, j \neq k(t)\big\}, 
$$
where the symbol $\downarrow$ means that we have ranked the terms in the decreasing order.
This defines a fragmentation process $H$ with parameters $(0,\nu)$, and all $(0,\nu)$-fragmentations can be constructed similarly from a Poisson point process. 

Next, from a homogeneous fragmentation $H$ with dislocation measure $\nu$ and any $a \in \mathbb R$, one can construct a fragmentation $F$ with parameters $(a,\nu)$ by using a family of suitable time-changes. And vice-versa, any  $(a,\nu)$-fragmentation can be turned into a homogeneous one with same dislocation measure by  suitable time-changes. We do not specify these time-changes here, as it requires some additional technical material, but refer to \cite{Ber02} for details. Let us just mention that the  time-changes  only depend on the past evolution of the masses. They modify the splitting times of the masses but neither change the masses themselves, nor  the genealogical relations between them. It turns out that for all $t \geq 0$ and all $i \in \mathbb N$, there exists a pair $(s,j) \in [0,\infty) \times \mathbb N$ such that $F_i(t)=H_j(s)$, and vice-versa.

\bigskip

\noindent \textbf{Fragmentation trees.} A \textit{rooted} continuous tree is a pair consisting of a continuous tree and one of its points, this distinguished point being called the \textit{root}. To simplify notation we write $ \mathrm{ht}(.)$ for the distance to the root in a rooted tree.

In \cite{HM04} it is proved that when $a<0$ and $\nu$ is conservative, if $F$ is an $(a,\nu)$ fragmentation process then one can construct a random compact rooted continuous tree $ \mathcal{T}^\bullet$ endowed with a (random) probability measure $\mu$ which is supported on its leaves such that $ \mathcal{T}^\bullet$ codes ``the genealogy'' of the process $F$ in the sense that almost surely for all $t \geq 0$ we have 
 \begin{eqnarray} F(t)&=&\downarrow \mu\text{-masses of the connected components of }
\{v \in \mathcal T^\bullet : \mathrm{ht}(v)>t\}. \label{etoile}
 \end{eqnarray}
The boundedness of this genealogical tree is due to a loss of mass by formation of dust in the fragmentation, which is itself due to the strictly negative index of self-similarity. The tree $ \mathcal{T}^\bullet$ is called the fragmentation tree associated with $F$. Recently, this result has been extended to the dissipative dislocation measures \cite{Steph}, the main difference being that the probability measure $\mu$ is then no more concentrated on the set of leaves of the genealogical tree and also charges its squeletton. When $\nu$ is dissipative and $\nu(0,0,...)=0$ it is even fully  supported by the squeletton.

A particular instance of fragmentation trees is given by the stable trees. Specifically, let $ \mathcal{T}_{\alpha}^\bullet$ denote an $\alpha$-stable tree endowed with its uniform mass measure $\mu_{\alpha}$ and rooted at $X_{0}$ (a uniform random leaf). If $F_{\alpha}$ denotes the $ \mathcal{S}$-valued process associated with $ \mathcal{T}^\bullet_{\alpha}$ by \eqref{etoile} 
then Bertoin \cite{Ber02} for the Brownian case and Miermont \cite{Mie03} for the general setting shew that $F_{\alpha}$ is a pure-jump self-similar fragmentation process with  self-similarity index $-\bar \alpha$ and computed its dislocation measure.
For the fragmentation $F_2$ this measure is given for all test functions $f:\mathcal S \rightarrow \mathbb R$ by
$$
\int_{\mathcal S}f(\mathbf s)\nu_{\mathrm{Br}}(\mathrm d \mathbf s)=\int_{1/2}^1 f(x,1-x,0,...) \left(\pi x^{3}(1-x)^{3} \right)^{-1/2}\mathrm dx,
$$
and for the stable fragmentation $F_{\alpha}$, $1<\alpha<2$, it is given by $$
\int_{\mathcal S}f(\mathbf s)\nu_{\alpha}(\mathrm d \mathbf s) =
C_{\alpha} \mathbb E\left[\sigma_1 f\left(\frac{\Xi_i}{\sigma_1},i \geq 1\right)\right],
$$
where 
\begin{equation}
\label{Cbeta}
C_{\alpha}=\frac{\alpha(\alpha-1) \Gamma(\bar \alpha)}{\Gamma(2 -
  \alpha)}
\end{equation}  
and $(\sigma_t, t\geq 0)$ is a stable subordinator of Laplace exponent 
$\lambda^{1/\alpha}$ and  $\left(\Xi_i,i\geq 1 \right)$  the sequence of
its jumps before time 1, ranked in the decreasing order.  Note that $\nu_{\alpha}(\mathcal S)=\infty$ for all $\alpha \in (1,2]$, which is related to the fact that the set of leaves of $\mathcal T_{\alpha}$ is dense in $\mathcal T_{\alpha}$. Note also that $\nu_{\alpha}$ is conservative.

We need some more notation. Note that for all $t\geq0$, the connected components  of $\{ v \in \mathcal{T}_{\alpha}^\bullet : \mathrm{ht}(v)>t\}$ are open subtrees of $\mathcal T_{\alpha}$.  We denote these subtrees by $ \mathcal{T}_{\alpha}^i(t), i \geq 1$ so that  $$ F_{\alpha}^i(t) = \mu_{\alpha}(  \mathcal{T}_{\alpha} ^i(t)),$$
with  $ F_{\alpha}(t)=(F_{\alpha}^i(t), i \geq 1)$. These subtrees are all distributed as $\mathcal T_{\alpha}\backslash \{X_0\}$, up to random scalings. Specifically, conditionally on $F_{\alpha}(t)$, the trees   $ \mathcal{T}_{\alpha}^i(t), i \geq 1$ are independent, such that for all $i$, $ \mathcal{T}_{\alpha}^i(t)$ endowed with the probability measure $\mu_{\alpha}(\ \cdot  \mid  \mathcal{T}_{\alpha}^i(t))$ is distributed as $F^i_{\alpha}(t) ^{\bar \alpha} \cdot (\mathcal T^i_{\alpha}\backslash\{X_0^i\})$ endowed with $\mu^i_{\alpha}$, where $(\mathcal T^i_{\alpha}, \mu^i_{\alpha})$ is an independent copy of the $\alpha$-stable tree and $X_0^i$ is a random point of  $\mathcal T^i_{\alpha}$ with distribution $\mu_{\alpha}^i$. This is called the self-similar property of the stable tree $\mathcal T_{\alpha}$. A similar property holds for all fragmentation trees.

\subsection{The tree $\mathfrak T_{\alpha,\alpha'}$ is a dissipative fragmentation tree}
\label{SecDiss}
From now on and in the remainder of this paper, we let $1<\alpha<\alpha' \leq 2$ and consider a stable tree $\mathcal T_{\alpha}$ together with a sample of uniform leaves $X_i,i \geq 0$.
We let $$\mathfrak{T}_{\alpha,\alpha'} = \mathrm{Prun}_{\alpha,\alpha'}( \mathcal{T}_{\alpha}, (X_{i})_{i\geq 0})$$ and denote by $ \mathcal{T}_{\alpha}^\bullet$ and $ \mathfrak{T}_{\alpha,\alpha'}^\bullet$ the trees  $ \mathcal{T}_{\alpha}$ and $ \mathfrak{T}_{ \alpha, \alpha'}$ rooted at the point $X_{0}$. Since $ \mathfrak{T}_{\alpha,\alpha'}$ is a closed subtree of $ \mathcal{T}_{\alpha}$ for every $x \in \mathcal{T}_{\alpha}$ there exists a unique point $ \mathrm{Proj}(x) \in \mathfrak{T}_{\alpha,\alpha'}$ minimizing the distance to $x$ in $ \mathfrak{T}_{\alpha,\alpha'}$. We then consider the measure $ \mu_{\alpha,\alpha'}$ on $ \mathfrak{T}_{\alpha,\alpha'}$ obtained as the push-forward of $\mu_{\alpha}$ by the application $ \mathrm{Proj}(.)$. 
\begin{proposition} The probability  $\mu_{\alpha,\alpha'}$ is purely atomic and supported by the skeleton of $ \mathfrak{T}_{\alpha,\alpha'}$.
\end{proposition}
\proof If $ x \notin \mathfrak{T}_{\alpha,\alpha'}$  then a moment of thought using the definition of the pruning procedure shows that $ \mathrm{Proj}(x)$ is a branch point of $ \mathfrak{T}_{\alpha,\alpha'}$. Thus the image of the measure $\mu_{\alpha}$ restricted to $ \mathcal{T}_{\alpha} \backslash \mathfrak{T}_{\alpha,\alpha'}$ is supported by the countable branch points of $ \mathcal{T}_{\alpha}$ that belong to $ \mathfrak{T}_{\alpha,\alpha'}$. The statement of the proposition follows since $ \mu_{\alpha}( \mathcal{T}_{\alpha} \backslash \mathfrak{T}_{\alpha,\alpha'}) =1$ by Proposition \ref{zeromass}.
\endproof 

The goal of this section is to show that $(\mathfrak T_{\alpha,\alpha'}^{\bullet},\mu_{\alpha,\alpha'})$ is the genealogical tree of a dissipative self-similar fragmentation process. To make this precise, we start by constructing its dislocation measure $\nu_{\alpha,\alpha'}$ from $\nu_{\alpha}$. For this, we design a procedure that extract randomly certain fragments of a dislocation of the unit mass.  

Let $ \mathbf{s} = (s_{1}, s_{2}, ...) \in \mathcal{S}$. We first consider $ \mathbf{s}^* = (s_{1}^*, s_{2}^*,...)$ a size-biased random reordering of this sequence from which we will extract a subsequence $ ( \widehat{s^*_{i}})_{i \geq 1}$ as follows. When $\alpha'=2$ the extraction is deterministic: set $\widehat{s^*_{1}}= s_{1}^* $ and $\widehat{s_{2}^*}= s_{2}^*$ and $\widehat{s_{n}^*}=0$ for $n\geq 3$. For $\alpha'<2$, we use an extra randomness and the definition of the probabilities $p_{\alpha,\alpha',d,d'}
$ defined by (\ref{pdd}) in the Introduction. We proceed recursively on $n$ as follows:
\begin{enumerate}
\item[(i)] $\widehat{s^*_{1}}= s_{1}^* $ and $\widehat{s_{2}^*}= s_{2}^*$,
\item[(ii)] for $n \geq 3$, let $d=n$ and $d'= \#\{ 1 \leq i \leq n-1 : \widehat{s_{i}^*}=s_{i}^*\}+1$ and put 
 \begin{eqnarray*}
\widehat{s_{n}^*}=s_{n}^* &\ \text{ with probability }& p_{\alpha,\alpha',d,d'},\\
\widehat{s_{n}^*}= 0 & \mbox{ otherwise}.
 \end{eqnarray*}
 \end{enumerate}
The dislocation measure $\nu_{\alpha,\alpha'}$ is then defined  for all test functions $f:\mathcal S
\rightarrow \mathbb R$ by
 \begin{eqnarray} \label{defnuaa}
\int_{\mathcal S
}f(\mathbf s)\nu_{\alpha,\alpha'}(\mathrm d \mathbf s)&: =& \int_{ \mathcal{S}} \nu_\alpha ( \mathrm{d} \mathbf{s}) \, \mathbb{E}\left[f\big( \downarrow( \widehat{s_{i}^*})\big)\right],
 \end{eqnarray}
where we recall that the symbol $ \downarrow$ means reordering in decreasing order. Note that $\nu_{\alpha,\alpha'}$ is  dissipative since $\alpha<\alpha'$. 
 
Our goal, now, is to prove the following result which clearly leads to Proposition \ref{thm:frag}.  The proof does not rely on the particular expression of the dislocation measure of the stable tree and the result could be generalized to any fragmentation tree pruned by a similar procedure.

\begin{proposition}
\label{propfrag} The process  $F_{\alpha,\alpha'}$  obtained by considering for each $t \geq 0$ the decreasing reordering of the 
\begin{eqnarray*}
&\mu_{\alpha,\alpha'}\text{-masses of the connected components of }
\big\{v \in \mathfrak T_{\alpha,\alpha'}^{\bullet} : \mathrm{ht}(v)>t\big\}
\end{eqnarray*}
is a self-similar fragmentation, with index $-\bar \alpha$ and dislocation measure $\nu_{\alpha,\alpha'}$.
\end{proposition}

To be coherent with the notation introduced for the fragmentation $ F_{\alpha}$ of the stable tree $\mathcal T_{\alpha}^{\bullet}$, we denote by $ \mathcal{T}_{\alpha, \alpha'}^i(t), i \geq 1$ the open subtrees forming the connected components of $\{ v \in \mathfrak{T}_{\alpha,\alpha'}^\bullet : \mathrm{ht}(v) >t\}$, so that the elements of $F_{\alpha,\alpha'}(t)$ satisfy for $i \geq 1$, $$ F_{ \alpha, \alpha'}^i(t) = \mu_{\alpha,\alpha'} ( \mathcal{T}_{\alpha,\alpha'} ^i(t)).$$ Note that for each $i \geq 1$, there exists a unique integer $j$ ($\geq i$) such that $ \mathcal{T}_{\alpha, \alpha'}^i(t) \subset  \mathcal{T}_{\alpha}^j(t)$, and then that
$$
 F_{ \alpha, \alpha'}^i(t) = \mu_{\alpha} ( \mathcal{T}_{\alpha} ^j(t))=F_{\alpha}^j(t).
$$ 
We then say that the mass $F_{\alpha}^j(t)$ contributes to  $F_{ \alpha, \alpha'}(t)$. Note that $F_{\alpha}^j(t)$ contributes to  $F_{ \alpha, \alpha'}(t)$ if and only if $\mathcal T_{\alpha}^j(t) \cap \mathfrak T_{\alpha,\alpha'} \neq \emptyset$.

\proof[Proof (Sketch)] A complete proof of the last proposition would be technical and we deliberately chose to stay rather informal. We first use the pruning construction of $\mathfrak T_{\alpha,\alpha'}$
to give a practical criterion to decide which of the elements of $F_{\alpha}(t)$ contribute to $F_{\alpha,\alpha'}(t)$. 
Let $F_{\alpha}^i(t)$ be an element of $F_{\alpha}(t)$ and consider the smallest $k \geq 1$ such that $X_{k} \in  \mathcal{T}_{\alpha}^i(t)$. Then $ F_{\alpha}^i (t)$ contributes to $F_{\alpha,\alpha'}(t)$ if and only if $X_{k} \in \mathrm{Prun}_{\alpha,\alpha'}( \mathcal{T}_{\alpha} ; X_{0}, ... , X_{k})$, that is if $X_k$ is a blue leaf, with the color coding of Section \ref{sec:pruning1}.
Consider then a jump of the process $F_{\alpha}$: Pick a time $t$ such that a fragment $F_{\alpha}^{i}(t-)$ splits into a countable collection of masses $$F_{\alpha}^{j}(t) \qquad  \mbox{ for } j \in J=\{ \mbox{indices offspring of }i \mbox{ at time }t\}.$$ In geometric terms, there exists a branch point at height $t$ in $ \mathcal{T}_{\alpha}^\bullet$ that branches into the subtrees $ \mathcal{T}_{\alpha}^j(t), j \in J$.  Let then assume that $F_{\alpha}^i(t-)$ contributes to $F_{\alpha,\alpha'}(t-)$. This means that during the pruning operation the first leaf $X_{k}$ falling into $ \mathcal{T}_{\alpha}^i(t-)$ is blue. 
Now, let us condition on that fact and on $ (F_{\alpha}(s), s \leq t)$ and set $$s_{j} = \frac{\mu_{\alpha}( \mathcal{T}_{\alpha}^j(t))}{\mu_{\alpha}( \mathcal{T}_{\alpha}^i(t-))}  \mbox{ for } j\in J.$$ Our goal is to show that the relative masses $s_k$ of the fragments $F_\alpha^k(t), k \in J$ that will contribute to  $ F_{\alpha,\alpha'}(t)$ are distributed as $(\widehat{s_j^*})_{j \in J}$.
We denote by $Y_{k}, k \geq 1$ the leaves among $X_{k}, k \geq 1$ that belong to $ \mathcal{T}_{\alpha} ^i(t-)$. Clearly, on our conditioning, the $(Y_{k})_{k \geq 1}$ are i.i.d. sampled according to $\mu_{\alpha}(\ \cdot \mid\mathcal{T}_{\alpha}^i(t-))$. Thus, if $ \mathcal{T}_{\alpha}^{a_{1}}(t), \mathcal{T}_{\alpha}^{a_{2}}(t), ...$ are the first, second, ... subtrees into which one of the points $Y_{i}$ falls then the sequence $s_{a_{1}}, s_{a_{2}}, ...$ is just a sized-biaised ordering of the sequence $(s_{j})_{j \in J}$, $$(s_{a_{n}})_{n\geq 1}  = (s_j^*)_{j\in J}.$$  Moreover, for every $n \geq 1$, the first leaf that falls into $  \mathcal{T}_{\alpha} ^{a_{n}}(t)$ is kept in the pruning construction (otherwise said is blue) with probability $p_{\alpha,\alpha',d,d'}$, where $d = a_{n}$ and $d'=\{ k \leq n-1 : \mathcal{T}_{\alpha}^{a_{k}}(t) \cap \mathfrak{T}_{\alpha,\alpha'} \ne \varnothing\} +1$. From this observation and the definition of $(\widehat{s_j}^*)$ it should be clear that the $s_{k}$'s, $k \in J$, corresponding to offspring subtrees of $\mathcal T^i_{\alpha}(t-)$ intersecting $ \mathfrak{T}_{\alpha,\alpha'}$ are thus distributed as $\widehat{s_{j}^*}, j \in J$ as desired. Obviously this procedure has to be done for \emph{each} jump of the fragmentation process. It is then possible to conclude by using the Poisson point process construction of the fragmentation $F_{\alpha}$. One important point is to notice that the time-changes used to pass to the homogenous counterparts of $F_{\alpha}$ and $F_{\alpha,\alpha'}$ are identical in both fragmentations for each element of $F_{\alpha}$ that contributes to $F_{\alpha,\alpha'}$. We leave subtleties to the careful reader.
\endproof

\subsection{From the dissipative fragmentation to a randomized conservative fragmentation}

Since the random tree  $\mathfrak{T}_{\alpha,\alpha'}$ is distributed as the stable tree $\mathcal T_{\alpha'}$ (up to a random scaling) it supports a uniform mass measure distributed as $\mu_{\alpha'}$.
However, the measure $\mu_{\alpha,\alpha'}$ (see last section) does not charge the set of leaves of the tree $ \mathfrak{T}_{\alpha,\alpha'}$ but only its skeleton and so is not connected, a priori,  to $\mu_{\alpha'}$.
We will see in this section how to construct from $\mu_{\alpha,\alpha'}$ a measure on the leaves of $\mathfrak{T}_{\alpha,\alpha'}^{\bullet}$ that will be related to $\mu_{\alpha'}$. 

\bigskip

\noindent \textbf{Keeping the mass.} Recall that the measure $\nu_\alpha$ is conservative, that is the mass is conserved at each splitting. However there is a loss of mass in $F_{\alpha}$ by formation of dust (due to the negative auto-similarity index). That is why, in this section, we rather focus on the homogeneous counterpart $H_\alpha$ of the fragmentation $F_{\alpha}$, which is obtained by applying suitable time-changes in $F_{\alpha}$, see Section \ref{sec:backgroundfrag}. In this case we have for all $t \geq 0$ $$ \sum_{i\geq1} H_{\alpha}^i(t) =1.$$
Obviously, this conservation of mass does not hold for homogeneous fragmentations with a dissipative dislocation measure. However there is then a natural way to define a random measure called the Malthusian measure which is stochastically conserved \cite{Ber06, BG04, Steph}. Although this new random measure is not a probability measure, it has a finite total mass with mean 1. 

Let us now introduce rigorously this object. For this we use  results established by Stephenson  \cite{Steph} for general dissipative fragmentation trees, following the ideas of Bertoin \cite{Ber06} for dissipative fragmentations with finite dislocation measures. From now on, we focus on the case of $ \nu_{\alpha,\alpha'}$ (see \eqref{defnuaa}), that is we consider $H_{\alpha,\alpha'}$ the homogeneous counterpart of the process $F_{\alpha,\alpha'}$. The starting point is the existence of a Malthusian index $p^*_{\alpha,\alpha'}$ which is characterized by
$$
\int_{\mathcal S} \bigg(1-\sum_{i \geq 1} s_i^{p^*_{\alpha,\alpha'}}\bigg) \nu_{\alpha,\alpha'}(\mathrm d \mathbf s)=0.
$$
See below for the discussion about existence and computation of this exponent. Once the existence of $p^*_{\alpha,\alpha}$ is granted, the process $$ \mathrm{Mart}_{1,0}(s) = \sum_{i\geq 1} \big(H^i_{ \alpha, \alpha'}(s)\big)^{p^*_{\alpha,\alpha'}}, \quad s\geq 0$$ is a positive martingale, see e.g. \cite{Ber06, Steph}. More generally, for each $t \geq 0$ and $i \in \mathbb{N}$, the process  \begin{eqnarray} \label{mart}\mathrm{Mart}_{i,t}(s) &=& \sum_{i \to j} \big( H^j_{ \alpha, \alpha'}(s)\big)^{p^*_{\alpha,\alpha'}}, \qquad s\geq t,  \end{eqnarray} where the sum is over all fragments whose ancestor at time $t$ is $H^i_{ \alpha, \alpha'}(t)$, is a positive càdlàg martingale which therefore converges almost surely towards a limiting value denoted by $ \mathrm{Mart}_{i,t}$, $i \in \mathbb{N}$, $t \geq 0$.  Actually, we may and will choose the random variables $\mathrm{Mart}_{i,t},i \in \mathbb N,t \geq 0$ so that the convergences hold  for \textit{all} $i,t$, almost surely (see  \cite{Steph}). Besides, by the homogeneity property of  the fragmentation process $H_{ \alpha, \alpha'}$, we clearly have for all $t\geq 0$, $$(\mathrm{Mart}_{i,t},i\geq 1) = \big((H^i_{ \alpha, \alpha'}(t) )^{p^*_{\alpha,\alpha'}}\mathrm{Mart}^{(i,t)}_{1,0},i \geq 0\big),$$ where the r.v. $\mathrm{Mart}^{(i,t)}_{1,0}, i\geq 0$ are independent, all distributed as $\mathrm{Mart}_{1,0}$ and also independent of $(H^i_{ \alpha, \alpha'}(t),i \geq 1)$.  Finally, recall the notation $\mathcal T_{\alpha}^i(t)$ and $\mathcal T_{\alpha,\alpha'}^i(t)$, $i \geq 1$ used in the last section for the open subtrees above level $t$ in $ \mathcal{T}_\alpha^\bullet$ and $ \mathfrak{T}_{ \alpha, \alpha'}^\bullet$ respectively. We similarly denote by $ \mathcal{H}_{\alpha}^{i}(t)$ and $ \mathcal{H}_{\alpha,\alpha'}^{i}(t)$ the open subtrees of $ \mathcal{T}_\alpha^\bullet$ and $ \mathfrak{T}_{ \alpha, \alpha'}^\bullet$  so that $$ H_\alpha^i(t)= \mu_\alpha( \mathcal{H}_\alpha^i(t)), \qquad H_{\alpha,\alpha'}^i(t)= \mu_{\alpha,\alpha'}( \mathcal{H}_{\alpha,\alpha'}^i(t)).$$
The set of subtrees $\{\mathcal H_{\alpha}^i(t), i \in \mathbb N,t \geq 0\}$ is in bijection with the set of subtrees $\{\mathcal T_{\alpha}^i(t), i \in \mathbb N,t \geq 0\}$, and it is the same for the subtrees indexed by $\alpha,\alpha'$.
We now have  the material to build the Malthusian measure $\mu^*_{ \alpha, \alpha'}$. On the event with probability 1 where all the martingales (\ref{mart}) converge,  we set 
 \begin{eqnarray}
\mu_{\alpha,\alpha'}^*\left(\mathcal H_{\alpha,\alpha'}^i(t)\right)&=&\mathrm{Mart}_{i,t}
 \label{def:mart-measure}
 \end{eqnarray}
for all $i \in \mathbb{N}$ and $t \geq0$.
By \cite{Steph}, this indeed defines uniquely a $\sigma$-finite measure $\mu^*_{\alpha,\alpha'}$ on the tree  $\mathfrak{T}_{\alpha,\alpha'}^{\bullet}$, which is fully supported by the set of leaves. Note that the set $\{\mathcal H_{\alpha,\alpha'}^i(t), i \in \mathbb N, t \geq 0\}$ generates the Borel $\sigma$-field on $\mathfrak T_{\alpha,\alpha'}$.

\bigskip

\noindent \textbf{Computing $p^*_{\alpha,\alpha'}$.} Apart from the Brownian case when $\alpha'=2$, it is not easy to calculate directly from the expression of $\nu_{\alpha,\alpha'}$ given in \eqref{defnuaa} its Malthusian index $p^*_{\alpha,\alpha'}$.  However we can use information from the proof of Proposition \ref{thm:inside}, based on the coupling in Marchal's construction, to get that
 \begin{eqnarray*}
p^*_{\alpha,\alpha'}&=&\frac{\bar \alpha
}{\bar \alpha
'}.
 \end{eqnarray*}
This will be done during the proof of Theorem \ref{thFrag}. There is another indirect way to get the value of $p^*_{\alpha,\alpha'}$, provided we know that this index exists and that the measure $\nu_{\alpha,\alpha'}$ satisfies suitable (weak) integrability properties. Indeed, on the one hand, the Hausdorff dimension of the $\alpha'$-stable tree, hence also that of  $\mathfrak{T}_{\alpha,\alpha'}^{\bullet}$, is $1/\bar \alpha
'$ a.s. \cite{DLG05,HM04}. On the other hand, according to results of \cite{Steph} on general dissipative fragmentation trees, we know that the Hausdorff dimension of the $(-\bar \alpha)$-self-similar fragmentation tree $\mathfrak{T}_{\alpha,\alpha'}^{\bullet}$ is almost surely $p_{\alpha,\alpha'}^*/\bar \alpha$, provided $\nu_{\alpha,\alpha'}$ satisfies the above mentioned  integrability properties. After identification this gives back $p^*_{\alpha,\alpha'}=\bar \alpha/\bar \alpha'$.

\bigskip 

The main result of this section is the fact that $\mathfrak{T}_{\alpha,\alpha'}^{\bullet}$ endowed with the measure $\mu^*_{\alpha,\alpha'}$ is distributed as an $\alpha'$-stable tree  endowed with its uniform mass measure, up to a random scaling.  Recall from (\ref{MomentsML}) the definition of a  generalized Mittag-Leffler random variable $\mathrm{ML}_{\bar \alpha/\bar \alpha',\bar \alpha}$ with parameters $({\bar{\alpha}}/\bar \alpha ',{\bar{\alpha}})$  and set $M_{\alpha,\alpha'}=(\alpha'/\alpha)^{1/\bar \alpha '}\mathrm{ML}_{\bar \alpha/\bar \alpha',\bar \alpha}$.
\begin{theorem} 
\label{thFrag}
There is the following identity in distribution, with $M_{\alpha,\alpha'}$ independent of $(\mathcal T_{\alpha'}^{\bullet}, \mu_{\alpha'})$,
 \begin{eqnarray*}
\Big(\mathfrak{T}_{\alpha,\alpha'}^{\bullet}, \mathbb E[M_{ \alpha, \alpha'}] \mu^*_{\alpha,\alpha'}\Big) &\overset{\mathrm{(d)}}=&\Big( M_{ \alpha, \alpha'}^{ \bar{ \alpha}'} \mathcal T_{\alpha'}^{\bullet}, M_{ \alpha, \alpha'} \mu_{\alpha'}\Big),
 \end{eqnarray*}
where $\mathbb E[M_{\alpha,\alpha'}]=(\alpha/\alpha')^{1/\bar \alpha
'} \bar \alpha
' \Gamma\big(\bar \alpha
\big)/\Gamma\big(\bar \alpha
(1+1/\bar \alpha
')\big)$, and where the measure $\mu_{\alpha '}$ in the right-hand side of this identity actually represents the push-forward of $\mu_{\alpha '}$ on $\mathcal T_{\alpha'}$ by the homothety that sends $\mathcal T_{\alpha'}$ on $M_{ \alpha, \alpha'}^{ \bar{ \alpha}'} \mathcal T_{\alpha'}$. 
\end{theorem}

In other words, the dissipative fragmentation tree  $\mathfrak{T}_{\alpha,\alpha'}^{\bullet}$,  endowed with a rescaled version of its Malthusian measure, is distributed as the genealogical tree of  an $\alpha'$-stable fragmentation of the random mass $M_{ \alpha, \alpha'}$. Note that the change of  measures  on the tree, from $\mu_{\alpha,\alpha'}$ to $\mu^*_{\alpha,\alpha'}$, changes the index of self-similarity of the underlying fragmentation, which passes from $-\bar \alpha$ to $-\bar \alpha
'$. \bigskip

\noindent \textit{Proof of Theorem \ref{thFrag}}. \textsc{Step 1.} Recall that $ \mathfrak{T}_{ \alpha, \alpha'} = \mathrm{Prun}_{ \alpha, \alpha'}( \mathcal{T}_{ \alpha}, (X_i)_{i\geq0})$ and that a leaf $X_i$ is said to be blue if it belongs to $ \mathfrak{T}_{ \alpha, \alpha'}$ and red otherwise. Write $ \mathbb{B}$ for the indices of the blues leaves and $ \mathbb{B}(n) = \mathbb{B} \cap \{0,1,2,...,n\}$. Recalling from the proof of Theorem \ref{thm:pruning} that $ \mathbb{B}(n)$ is the number of leaves of  a chain $\mathbf T_{\alpha,\alpha'}$, we apply  Lemma \ref{TnRRemi} to get 
 \begin{eqnarray}
\label{cvM}
\frac{\# \mathbb{B}(n)}{n^{\bar \alpha
/\bar \alpha
'}} &\xrightarrow[n\to\infty]{\mathrm{a.s.}}& \mathrm{ML},
 \end{eqnarray}
where $\mathrm{ML}$ has a generalized Mittag-Leffler distribution with parameters $(\bar \alpha
/\bar \alpha',\bar \alpha)$. In particular, its mean is 
$$
\mathbb E[\mathrm{ML}]=\frac{\bar \alpha
' \Gamma\left(\bar \alpha
\right)}{\Gamma\big(\bar \alpha
\big(1+1/\bar \alpha
'\big)\big)},
$$
using (\ref{MomentsML}) and the relation $\Gamma(a+1)=a\Gamma(a)$ to get simplifications. Moreover, we know from (the proof of) Proposition \ref{thm:inside} and the discussion at the end of Section \ref{secproofpruning} that there exists an $\alpha'$-stable tree $\mathcal T_{\alpha'}$ such that
$$ \mathfrak{T}_{\alpha, \alpha'} =\frac{\alpha'}{\alpha} \mathrm{ML}^{\bar \alpha '}\mathcal T_{\alpha'}$$ 
and that
$(X_i)_{i \in \mathbb{B}}$ is an i.i.d. sample of the uniform mass measure on (the leaves of) $\mathfrak{T}_{\alpha, \alpha'}$, which, with a slight abuse of notation, is also denoted by $\mu_{\alpha'}$.
Consequently, 
\begin{equation}
\label{empirique}
 \mu_{\alpha'} \text{ 
is the almost sure empirical limit of} \quad \frac{\sum_{i \in \mathbb{B}(n)} \delta_{X_i}}{\#  \mathbb{B}(n)} \quad \text{as }n \rightarrow \infty,
\end{equation}
which will be useful in the following. Theorem \ref{thFrag} will therefore be proved if we show the existence of $p^*_{\alpha,\alpha'}$ (hence that of $\mu^*_{\alpha,\alpha'}$) and that almost surely for all $i \in \mathbb N, t \geq 0$,
\begin{equation}
\label{eqidentity}
 \mathrm{ML} \mu_{\alpha'} \left(\mathcal H_{\alpha,\alpha'}^i(t) \right)= \mathbb E[\mathrm{ML}]\mu_{\alpha,\alpha'}^*\left(\mathcal H_{\alpha,\alpha'}^i(t)\right).
\end{equation}
For this, we will use (\ref{cvM}) and (\ref{empirique}).
\medskip

\noindent \textsc{Step 2.} As their self-similar counterparts, $ H_{\alpha,\alpha'}$ is a sub-fragmentation of $H_{\alpha}$.
More precisely, for $i \in \mathbb{N}$ and $t \geq 0$, there exists a unique integer $j$ such that $ \mathcal H_{\alpha,\alpha'}^i(t) \subset \mathcal H_{\alpha}^j(t)$ and moreover, 
$$
H^i_{\alpha,\alpha'}(t)=\mu_{\alpha,\alpha'}( \mathcal H_{\alpha,\alpha'}^i(t))=\mu_{\alpha}( \mathcal H_{\alpha}^j(t) )=H^j_{\alpha}(t).
$$
Then, similarly to (\ref{cvM}), introduce  
the random variable 
\begin{eqnarray*}
\mathrm{ML}_i(t)&=&\lim_{n \rightarrow \infty} \frac{\# \{ k \in \mathbb{B}(n) : X_{k} \in \mathcal{H}_{\alpha}^{j}(t)\}}{\# \{ k \leq n : X_{k} \in \mathcal{H}_{\alpha}^{j}(t)\} ^{ \bar{ \alpha}/ \bar \alpha'}}
\end{eqnarray*}
which, by the self-similarity property of  $\mathcal T_{\alpha}$ and of the blue-coloring operation exists a.s., is distributed as $\mathrm{ML}$ and is independent of $H_{\alpha}(t)$. Moreover, the random variables $\mathrm{ML}_{i}(t), i \geq 1$ are independent. Now note the following easy but crucial fact: almost surely,
\begin{eqnarray}
\label{eg1}
\mathrm{ML}  \mu_{\alpha'} \left(\mathcal H_{\alpha,\alpha'}^i(t) \right) &=&\lim_{n \rightarrow \infty} \frac{\#  \mathbb{B}(n)}{n^{\bar \alpha
/\bar \alpha
'}} \times  \frac{\# \{ k \in \mathbb{B}(n) : X_{k} \in \mathcal{H}_{\alpha}^{j}(t)\} }{\#  \mathbb{B}(n)} \nonumber \\
&=& \lim_{n \rightarrow \infty} \frac{\# \{ k \in \mathbb{B}(n) : X_{k} \in \mathcal{H}_{\alpha}^{j}(t)\}}{\# \{ k \leq n : X_{k} \in \mathcal{H}_{\alpha}^{j}(t)\} ^{ \bar{ \alpha}/ \bar \alpha'}}\times \left(\frac{\{ k \leq n : X_{k} \in \mathcal{H}_{\alpha}^{j}(t)\}}{n}\right)^{\bar \alpha/\bar \alpha'}\nonumber \\
&=& \mathrm{ML}_i(t) \left(H^{j}_{\alpha}(t)\right)^{\bar \alpha
/\bar \alpha
'} = \mathrm{ML}_i(t) \left(	H^i_{\alpha,\alpha'}(t)\right)^{\bar \alpha/\bar \alpha'}.\end{eqnarray}
Note also that a.s.
$$
\sum_{i \geq 1}\mu_{\alpha'}\left(\mathcal H_{\alpha,\alpha'}^i(t) \right)=1, \quad \forall t \geq 0.$$ 
Indeed, since $H_\alpha$ is the homogeneous version of the fragmentation $F_\alpha$, the union $\cup_i \mathcal{H}_\alpha ^i(t)$ contains all the leaves of $ \mathcal{T}_\alpha$ but the root, $\forall t \geq 0$ a.s. Consequently, $\cup_i \mathcal{H}_{\alpha,\alpha}^i(t)$ also contains all the leaves (except the root) of $ \mathfrak{T}^\bullet_{ \alpha, \alpha'}$, hence the last display. Together with (\ref{eg1}), this leads to
\begin{equation}
\label{eqfondamentale}
\mathrm{ML}=\sum_{i \geq 1} \mathrm{ML}_i(t) \left( H_{\alpha,\alpha'}^i(t)\right)^{\bar \alpha
/\bar \alpha
'}, \quad \forall t\geq 0 \text{ a.s.}
\end{equation} 
Taking expectations and recalling that $0<\mathbb E[\mathrm{ML}]<\infty$, this implies that
$$
1=\mathbb E \bigg[\sum_{i \geq 1}  \left( H_{\alpha,\alpha'}^i(t)\right)^{\bar \alpha
/\bar \alpha'}\bigg]= \exp\bigg(- t\int_{\mathcal S
}\bigg(1-\sum_{i\geq 1} s_i^{\bar\alpha/\bar \alpha'} \bigg) \nu_{\alpha,\alpha'}
(\mathrm d \mathbf s)\bigg),
$$
using for the second equality that $\nu_{\alpha,\alpha'}$ is the dislocation measure of $H_{\alpha,\alpha'}$ and then the well-known fact in the theory of homogeneous fragmentations
$$
\mathbb E \bigg[\sum_{i \geq 1}  \left( H_{\alpha,\alpha'}^i(t)\right)^{q}\bigg]= \exp\bigg(- t\int_{\mathcal S
}\bigg(1-\sum_{i\geq 1} s_i^{q}\bigg) \nu_{\alpha,\alpha'}
(\mathrm d \mathbf s)\bigg), \quad \forall q \in \mathbb R, 
$$
see e.g. Bertoin \cite{Ber01}. Thus we indeed have $p^*_{\alpha,\alpha'}=\bar \alpha
/\bar \alpha'$, and then the existence of the almost sure limits  $\mathrm{Mart}_{i,t}$ of the martingales (\ref{mart}), and also that of the measure $\mu^{*}_{\alpha,\alpha'}$.

\medskip

\noindent \textsc{Step 3.}  It remains to use (\ref{eg1}) to get (\ref{eqidentity}). Following an idea of \cite{BG04} we get from (\ref{eqfondamentale}) and the definition of the martingale $(\mathrm{Mart_{1,0}}(t),t \geq 0)$ that
\begin{eqnarray*}
\mathbb E \left[\big(\mathrm{ML}-\mathbb E[\mathrm{ML}]\mathrm{Mart_{1,0}}(t) \big)^2 \right] &=& 
\mathbb E\Bigg[ \bigg(\sum_{i\geq 1}  \left( H_{\alpha,\alpha'}^i(t)\right)^{\bar \alpha
/\bar \alpha
'}\left(\mathrm{ML}_i(t)-\mathbb E[\mathrm{ML}] \right)\bigg)^2\Bigg] \\
&=&\mathrm{Var}(\mathrm{ML})\mathbb E \Bigg[ \sum_{i \geq 1}  \left( H_{\alpha,\alpha'}^i(t)\right)^{2\bar \alpha
/\bar \alpha
'} \Bigg] \\
&=& \mathrm{Var}(\mathrm{ML}) \exp\Bigg(- t\int_{\mathcal S
}\Bigg(1-\sum_{i\geq 1} s_i^{2\bar \alpha
/\bar \alpha
'} \Bigg) \nu_{\alpha,\alpha'}
(\mathrm d \mathbf s)\Bigg),
\end{eqnarray*}
where the second equality comes from the fact that the random variables $\mathrm{ML}_i(t),i \geq 1$ are i.i.d distributed as $\mathrm{ML}$, and independent of $H_{\alpha}(t)$. The function 
$$q \mapsto \int_{\mathcal S
}\Bigg(1-\sum_{i\geq 1} s_i^q \Bigg) \nu_{\alpha,\alpha'}
(\mathrm d \mathbf s)$$
is strictly increasing on its domain of definition and is equal to 0 for $q=\bar \alpha
/\bar \alpha
'$, as seen above. Hence, noticing that $\mathrm{Var}(\mathrm{ML})<\infty$,
$$\mathbb E \left[\big(\mathrm{ML}-\mathbb E[\mathrm{ML}]\mathrm{Mart_{1,0}}(t) \big)^2 \right] \rightarrow 0\quad \text{as }t \rightarrow \infty,$$
which implies the almost sure equality
$$
\mathrm{ML}=\mathbb E[\mathrm{ML}]  \mathrm{Mart}_{1,0}, 
$$
since moreover $\mathrm{Mart_{1,0}}(t) \rightarrow \mathrm{Mart}_{1,0}$ a.s.
Similarly, by the homogeneity property of the fragmentation $H_{\alpha}$, for each $(i,t)$,
$$
\mathrm{ML}_{i}(t)\left( H_{\alpha,\alpha'}^i(t)\right)^{\bar \alpha
/\bar \alpha
'}=\mathbb E[\mathrm{ML}]  \mathrm{Mart}_{i,t}, \quad \text{a.s.},
$$ 
which, together with (\ref{eg1}), finally gives the expected
$$
\mathrm{ML} \mu_{\alpha'}\left(\mathcal H_{\alpha,\alpha'}^i(t) \right) =\mathbb E[\mathrm{ML}]  \mathrm{Mart}_{i,t}=\mathbb E[\mathrm{ML}]  \mu_{\alpha,\alpha'}^*\left(\mathcal H_{\alpha,\alpha'}^i(t)\right) \quad \text{a.s.}
$$
Note that this identity is proved for each fixed $(i,t)$ almost surely, and we actually want it almost surely for all $(i,t) \in \mathbb N \times [0,\infty)$. This is indeed true: for each  $(i,t)$, consider a sequence $(t_n)$ with $t_n \in \mathbb Q \cap [0, \infty)$ decreasing towards $t$; then one has 
$$
\mathcal H_{\alpha,\alpha'}^i(t)= \bigcup_{n \geq 1} \bigcup_{j \in J_n(i,t)}  \mathcal H_{\alpha,\alpha'}^{j}(t_n)
$$
where $J_n(i,t)$ is the set of all indices $j$ such that the mass  $H_{\alpha,\alpha'}^j(t_n)$ is a descendant of $H_{\alpha,\alpha'}^i(t)$. The above union on $j \in J_n(i,t)$ is disjoint, whereas that on $n$ is an increasing union. Consequently, a measure on $\mathfrak{T}_{\alpha,\alpha'}$ is entirely characterized by the weights it assigns to the subtrees $\mathcal H_{\alpha,\alpha'}^j(q), j \in \mathbb N, q \in \mathbb Q \cap [0,\infty)$. This ends the proof (and the paper). $\hfill \square$

\bibliographystyle{siam}

\end{document}